\documentclass[10pt]{article}
\usepackage{latexsym,amssymb}
\usepackage{t1enc}
\usepackage[cp1250]{inputenc}
\def\d{\,\mathrm{d}}
\def÷{\mathrm{div}}

\def\ind{\mathrm{ind}}
\def\Ind{\mathrm{ind}}
\def\deg{\mathrm{deg}}
\def\Deg{\mathrm{Deg}}

\def\wew{\mathrm{int}\,}

\def\part{\partial}
\def\bd{\partial}

\def\Ker{\mathrm{Ker}}

\def\R{\mathbb{R}}

\def\Z{\mathbb{Z}}

  \newcommand{\lma}{\lambda}
  
  \newcommand{\eps}{\varepsilon}
  \newcommand{\be}{\begin{equation}}
  \newcommand{\ee}{\end{equation}}
  \newcommand{\bp}{ \begin{pkt} \em }
  \newcommand{\ep}{ \end{pkt} }

\begin{document}

\thispagestyle{empty}

\begin{center}

$\mbox{ }$

\vspace{40mm}

\begin{Large} \begin{sc} 
Periodic and stationary solutions of nonlinear evolution equations -- translation along trajectories method
\end{sc}\end{Large}

\vspace{10mm}

{\bf Self-report }\\

\vspace{40mm}

{\bf Aleksander \'{C}wiszewski}

\vspace{5mm}

\noindent  Nicolaus Copernicus University in Toru\'{n}\\
Faculty of Mathematics and Computer Science\\
ul. Chopina 12/18, 87-100, Toru\'{n}\\
e-mail: Aleksander.Cwiszewski@mat.umk.pl

\vspace{50mm}

Toru\'{n}, October 2011\\

\end{center}

\vspace{0mm}

\newpage

$\mbox{ }$

\thispagestyle{empty}

\newpage

$\mbox{ }$

\vspace{30mm}

\noindent  {\bf Contents}

\noindent 1. Introduction 

\noindent 2. Compactness in evolution equations % Zwartość w równaniach ewolucyjnych

\noindent 3. Krasnosel'skii type theorems % Twierdzenia typu Krasnosielskiego

\noindent 4. Averaging principle for periodic solutions % Zasada uśredniania dla rozwiązań okresowych

\noindent 5. Averaging principle for hyperbolic equations at resonance % Zasada uśredniania dla równań hiperbolicznych w rezonansie

\noindent 6. Poincar\'{e}-Hopf theorem in infinite dimension % Twierdzenie Poincar\'{e}go-Hopfa w nieskończonym wymiarze

\noindent 7. Averaging and periodic solutions for partial differential equations

\indent a. Periodic solutions for damped nonlinear hyperbolic equations

\indent b. Periodic solutions of strongly damped beam equation

\indent c. Nonnegative periodic solutions for parabolic problems

\noindent 8. Translation along trajectories method vs. other methods of studying periodic problems

\noindent 9. Summary of results not included in the habilitation dissertation % Omówienie wyników niewchodzących w skład rozprawy

\noindent 10. About the Author

References

% \vspace{10mm}

$\mbox{ }$

\vspace{15mm}

\newpage

$\mbox{ }$

\thispagestyle{empty}

\newpage

$\mbox{ }$

\vspace{15mm}

%\vspace{-1.5mm}

{\bf 1. Introduction}

This work is a presentation of the results that constitute my habilitation dissertation.
It is devoted to applications of topological methods to nonlinear evolution equations.
The main concept of the approach proposed here is Poincar\'{e}-Krasnosel'skii operator of translation along trajectories, which is used to examine the existence and behavior of both periodic and stationary solutions of differential equations. 
Nonlinear evolution equations constitute an abstract setting -- in terms of functional analysis -- for many types of equations and systems of partial differential equations: many equations evolving from physics, chemistry and various mathematical models of other fields in science and technology, among others.
Therefore, evolutions equations are subject of intensive studies.\\
\indent Topological methods gave a series of effective results in the theory of ordinary differential equations. One of them is the concept of translation along trajectories operator -- also known as Poincar\'{e}-Krasnosel'skii operator. It goes back as far as H. Poincar\'{e} and was significantly developed by M.A. Krasnosel'skii.
Using topology to investigate the properties of translation operator leads to nontrivial quantitative and qualitative results relating existence of equilibrium points, periodic solutions as well as their bifurcations and branching.\\
\indent The goal of my dissertation is to examine properties of the translation operator for evolution equations in Banach spaces and to use it to studying periodic solutions and  stationary solutions.
However, in the case, the infinite dimensional nature and the related lack of compactness require using sophisticated concepts of functional analysis and a suitable development of the theory of homotopy invariants.

The dissertation consists of the following publications:

\noindent {\bf [R1]} {\em Krasnosel'skii type formula and
translation along trajectories method for evolution equations},
Discrete and Continuous Dynamical Systems Ser. A, Vol. 22, No. 3
(2008), 605--628. Co-author: P. Kokocki;

%\vspace{-2mm}

\noindent {\bf [R2]} {\em Periodic solutions for nonlinear
hyperbolic evolution systems}, Journal of Evolution Equations 10
(2010), 677--710. Co-author: P. Kokocki;

%\vspace{-2mm}

\noindent {\bf [R3]} {\em On a generalized Poincar\'{e}-Hopf
formula in infinite dimensions}, Discrete and Continuous Dynamical
Systems Series A, vol. 29 no. 3 (2011), 953--978. Co-author: W.
Kryszewski;

%\vspace{-2mm}

\noindent {\bf [R4]} {\em  Positive periodic solutions of parabolic evolution
problems: a translation along trajectories approach}, Central European Journal of
Mathematics, vol. 9, no. 2 (2011), 244--268;

%\vspace{-2mm}

\noindent {\bf [R5]} {\em Periodic solutions of damped hyperbolic
equations at resonance: translation along trajectories approach},
Differential an Integral Equations, vol. 24, no. 7-8 (2011),
767--786;

%\vspace{-2mm}

\noindent{\bf  [R6]} {\em Periodic oscillations for strongly
damped hyperbolic beam equation}, Topological Methods in Nonlinear
Analysis, vol. 37, no. 2, (2011), 259--282.

% \vspace{-3mm}

In my opinion the most significant results of the dissertation are the following:

\vspace{-1mm}

\begin{itemize}
\item the Krasnosel'skii type formulae for nonlinear perturbations of $C_0$ semigroup generators: [R1, Th. 5.1, Th. 5.3], [R4, Th. 4.1] and [R6, Th. 3.1];
\item averaging and continuation principles for periodic solutions of general nonlinear evolution problems: 
[R1,Prop. 4, Th. 6.2], [R2, Th 4.4, 4.10], [R4, Th. 5.5, Th. 5.7] and [R6, Th. 4.2];
\item the averaging and continuation principles for damped hyperbolic problems at resonance: [R5, Th. 3.1, Th. 3.7];
Landesman-Lazer type criteria for the existence of periodic solutions for damped hyperbolic partial differential equations: [R5, Th. 4.1];
\item the infinite dimensional Poincar\'{e}-Hopf type formula for abstract evolution equations 
determined by nonlinear perturbations of generators of compact $C_0$ semigroups: [R3, Th. 1.2, Th. 5.2].
\end{itemize}

In the remainder of Section 1 I present a sketch of ideas and general statements of results as well as
I point out the main difficulties that come up in their proofs.
One of most important matters, related with topological methods in evolution equations, 
is the compactness of translation along trajectories operator.
Some types of compactness, that are common in partial differential equations, are illustrated in Section 2 
with concrete relevant examples. The method of finding periodic solutions proposed in the dissertation
is based on a Krasnosel'skii type formula, that I present in Section 3, and an averaging principle, discussed in Section 4. The developed methods are applied to abstract hyperbolic equations of second order at resonance, for which 
a proper continuation method is provided in Section 5.
The obtained results enabled us to establish the relation between topological degree and Conley homotopy index, to which I devote Section 6. The presented general results for evolution equations are employed for partial differential equations.
The criteria for the existence of periodic solutions are given in Section 7.
Section 8 sheds some light on pros and cons of translation along trajectories and other approaches used for periodic problems in partial differential equations. Finally, in Section 9, I make a brief presentation of my results that 
are not included in the habilitation dissertation. A few facts about the author are enlisted in Section 10.

We are concerned with differential equations of the form
$$
\dot u (t) = A u(t) + F(t,u(t)),\ \ t>0         \leqno{(Z)} 
$$ where the linear operator $\!A\!:\!D(A)\!\to\!\mathbb{E}\!$ is the generator 
of a $C_0$ semigroup $\{S_A(t)\!:\!\mathbb{E}\!\to\! \mathbb{E} \}_{t\geq 0}$ 
of bounded linear operators on a Banach space $\mathbb{E}$ and the perturbation
$F:[0,+\infty)\times \mathbb{E} \to \mathbb{E}$ is continuous and $T$-periodic with respect to the first variable.
It is noteworthy at this point that many equations and systems of partial differential equations have this form
after a reformulation. It is the case for parabolic equations (in particular, for reaction diffusion systems and heat equations) and hyperbolic equations (in particular, for telegraph, beam, string equations and their versions with damping and strong damping). The main tool in studying problems of the form $(Z)$ is the well developed theory of $C_0$ semigroups (see \cite{Pazy}, \cite{Tanabe}), which, under additional assumptions on $F$,  assures the existence and uniqueness properties for the equation $(Z)$.
It is well known that under some additional conditions on $F$,
problem $(Z)$ admits the existence and uniqueness properties for solutions satisfying initial value condition.
One considers here  so-called {\em mild solutions}, i.e. continuous functions $u:[0,+\infty)\to \mathbb{E}$ 
satisfying
\be\label{15052011-2004} u(t) = S_A(t) \bar
u(0) + \int_{0}^{t} S_A  (t-s) F(s,u(s)) \,\d s \ \ \mbox{ for all } t>0.
\ee
We do not discuss here a completely separate issue of regularity of solutions and in the sequel
mild solutions will be referred to as {\em solutions}.
If we assume that $F$ satisfies the local Lipschitz condition, the equation $(Z)$ admits
the local existence and uniqueness properties, and under some growth conditions, dissipativity of the equation or the existence of a proper potential, $(Z)$ admits the global existence (on the half line $[0,+\infty)$) and uniqueness properties.
Throughout the whole self-report we assume that equations possess the global existence and uniqueness properties.
Then, we can correctly associate with $(Z)$ a mapping  
$\Phi_t:\mathbb{E}\to \mathbb{E}$, $t>0$, given by 
$\Phi_t (\bar u):= u(t)$, $\bar u\in \mathbb{E}$, where $u$ is the solution of $(Z)$ satisfying the condition $u(0)=\bar
u$. This mapping is called the {\em translation along trajectories operator} by time $t$.
It is worth pointing out that for the autonomous equation (if $F$ does not depend on time $t$), the family of mappings
$\{ \Phi_t\}_{t\geq 0}$ determines a {\em semiflow}
({\em semidynamical system } or {\em semigroup}) on the space $\mathbb{E}$.\\
\indent Clearly, the existence of a fixed point $\bar u$ for $\Phi_T$, $T>0$, i.e. the equality
$\Phi_T (\bar u) = \bar u$, implies the existence of a $T$-periodic solution for the problem $(Z)$.
Hence, to get $T$-periodic solutions one needs criteria for the existence of fixed points for the operator $\Phi_T$. 
The key issue is the compactness property of $\Phi_T$, since it is necessary in topological fixed point theorems (either in global or local versions). It appears that, under additional compactness assumptions on the linear operator $A$ and the perturbation $F$, some kind of their compactness is inherited by $\Phi_T$.
Then one can apply appropriate fixed point theorem, e.g. Schauder's Theorem -- when $\Phi_T$ is completely continuous,
Sadovskii's Theorem -- when $\Phi_T$ is condensing with respect to a measure of noncompactness or
Leschetz's Theorem -- if we search for solutions in a constraint set $M\subset \mathbb{E}$ being a~neighborhood retract 
in $\mathbb{E}$. This approach gave rise to numerous results on periodic solutions, see for instance the papers by Browder \cite{Browder1965}, Becker \cite{Becker}, Pr\"{u}ss \cite{Pruss}, Dancer and Hess \cite{Dancer-Hess} and \cite{Dancer-Hess1}, Hess \cite{Hess},
Shioji \cite{Shioji1997},  Bothe \cite{Bothe} and \cite{Bothe1}, Bader and~Kryszewski \cite{Bader-Kryszewski}, as well as the author \cite{JDE2009-Cw-Kr}.\\
\indent In the papers included in the dissertation, instead of global fixed point theorems there are used more subtle tools such as {\em local homotopy invariant}, precisely: {\em fixed point index theory} adequate to the type of a problem under consideration. The nontriviality of the fixed point index $\ind (\Phi_T, U)$, with respect to a subset $U \subset \mathbb{E}$, implies the existence of a fixed point in $U$ and, in consequence, a periodic solution of the problem $(Z)$.
Therefore, when considering the periodic problem, the key point is to determine the fixed point index of 
the translation operator $\Phi_T$.
Krasnosel'skii type theorems and a version of averaging principle will play important roles.
They allow us to express the fixed point index of the translation operators by means of an adequate homotopy index for the right hand side of the equation.\\
\indent The classical Krasnosel'skii theorem deals with an ordinary differential equation
$\dot u(t)  = f(u(t))$, $t\geq 0$, with $f:\R^N\to \R^N$ being a continuous mapping such that 
the existence and uniqueness properties hold. It states that, for an open bounded set $U\subset \R^N$ with no zeros of $f$ in its boundary, the index of the translation operator $\Phi_t$ with respect to $U$ is equal to the Brouwer degree $\deg_B(-f, U)$ for sufficiently small $t>0$ (see \cite{Krasnosielski}, \cite{Krasnosielski-Zabreiko}).\\ 
\indent In the dissertation I~investigate an autonomous equation
\be\label{15042011-1914} \dot u (t) = A u(t)
+ F(u(t)),\ \ t>0, 
\ee 
where, as before, $A$ is a $C_0$ semigroup generator and
$F:\mathbb{E}\to \mathbb{E}$ is a continuous perturbation.
Under additional compactness assumptions, I proved that given an open bounded set
$U\subset \mathbb{E}$ such that $A\bar u + F(\bar u) \neq
0$ for $\bar u \in
\partial U \cap D(A)$, there exists $t_0>0$ such that for all $t\in (0, t_0]$, $\Phi_t (\bar u)\neq \bar u$, for any $\bar u\in
\partial U$, and
$$
\Ind (\Phi_t, U) = \Deg ( A+F, U)
$$
where $\Deg(A+F, U)$ is the proper version of topological degree (see Remark 3.2).
The above formula was obtained in a few versions, under various compactness conditions -- adequate to studied 
equations. The results are contained in [R1], [R4] as well as [R6].\\
\indent Let us return to the equation $(Z)$. In order to determine $\ind (\Phi_T, U)$ we apply a~proper averaging principle. Let $u(\cdot; \bar u, \lambda):[0,+\infty)\to \mathbb{E}$, for $\bar u\in \mathbb{E}$ and a parameter $\lma>0$, stands for the solution of a problem
$$
\dot u(t) = A u(t) + F(t/\lma, u(t)), \ t>0, \leqno{(Z_\lma)}
$$
with the initial value condition $u(0) = \bar u$. And let $\widehat u (\cdot;\bar u):[0,+\infty)\to \mathbb{E}$ 
be the solution of the averaged autonomous equation
$$
\dot u(t) = A u(t) + \widehat F(u(t)),     \ t>0, \leqno{(\widehat
Z\, )}
$$
with the initial value condition $u (0)=\bar u$, where 
the mapping $\widehat F:\mathbb{E}\to \mathbb{E}$ is given by 
$$
\widehat F (\bar u):= \frac{1}{T} \int_{0}^{T}F(t,\bar u) \d t, \, \bar u\in \mathbb{E}.
$$
The averaging principle states that if the mapping $F$ is $T$-periodic in the first variable, then
$$
u(t; \bar v, \lma ) \to \widehat u (t;\bar u) \mbox{ as  } \bar v\to \bar u \mbox{ and } \lma \to 0,
$$
uniformly with respect to $t$ from bounded intervals.
I proved that principle in [R4] for general $C_0$ semigroup generators and perturbations $F$ satisfying a general compactness condition. The theorem is an infinite dimensional generalization 
of classical theorems by Bogolyubov and Mitropolsky (\cite{Bogo-Mitro})
as well as Henry's results (\cite{Henry}) and the work by Couchouron and Kamenskii (\cite{Couchouron-Kamenski}) (for more details see the comment in Section 4, page \pageref{15062011-1001}).\\
\indent Now pass to the periodic problem associated with $(Z)$, that is
$$
\left\{
\begin{array}{l}
\dot u (t) = A u(t) + F(t, u(t)),\ t>0,\\
u(0)=u(T),
\end{array}
\right.
$$
where the mapping $F:[0,+\infty)\times \mathbb{E}\to \mathbb{E}$ is $T$-periodic in the first variable.
In order to find the index $\ind (\Phi_T, U)$, I employ here a scheme that was used earlier for ordinary differential equations, e.g. \cite{Furi-Pera} and \cite{Furi-Pera-Spadini}.
Its main idea is that we embed the problem $(Z)$ into the family of problems 
$(Z_\lma)$, $\lma\in (0,1]$, next we use averaging methods for small $\lma\in (0,1]$ 
and continuation to get the existence of solutions for $\lma=1$. I assume that for an open bounded 
$U\subset\mathbb{E}$
$$
A\bar u + \widehat F (\bar u) \neq 0 \mbox{ for } \bar u\in \partial U\cap D(A)
$$
and show that, for sufficiently small $\lma>0$, the translation operator $\Phi_{\lma T}^{(\lma)}$ 
associated with  the equation $(Z_\lma)$ is homotopic to the translation operator $\widehat \Phi_{\lma T}$ of the  averaged equation $(\widehat Z\,)$. In consequence, due to the  homotopy invariance of fixed point index
\be\label{15062011-1016}
\ind (\Phi_{\lma T}^{(\lma)}, U) = \ind (\widehat \Phi_{\lma T}, U) \mbox{ for sufficiently small } \lma>0.
\ee
Next, applying the mentioned Krasnosel'skii type formula to $(\widehat Z)$, one obtains
$$
\ind (\widehat \Phi_{\lma T}, U) = \Deg ((A,\widehat F), U) \mbox{ for small } \lma>0,
$$
which, combined with (\ref{15062011-1016}), yields a formula
$$
\ind (\Phi_{\lma T}^{(\lma)}, U) = \Deg ((A,\widehat F), U) \mbox{ for small } \lma>0.
$$
In the case when proper {\em a priori } type estimates are satisfied, the homotopy invariance of the index implies
that  $\ind(\Phi_T, U) = \Deg((A,\widehat F), U)$ (as $\Phi_{T}^{(1)} = \Phi_T$).
In this manner, it is proved that periodic solutions of $(Z_\lma)$ emanate from equilibria of the equation $(\widehat Z\,)$ having nontrivial topological index and {\em a priori } type estimates enable us to ''track''
(i.e. localize) branches of these solutions and, therefore, get $T$-periodic solution of the original problem $(Z)$.
The already sketched approach is employed in [R1], [R2], [R4] and [R6].\\
\indent In the situation when the equation $(Z)$ is at resonance, i.e. $\Ker A \neq \{ 0 \}$ and $F$ is a bounded continuous perturbation, the above procedure does not apply. 
Then we consider a family of problems with parameters  $\eps>0$ and $\lma>0$
$$
\dot u (t) = Au(t) + \eps F(t/\lma, u(t)), t>0. % \leqno{(Z_{\lma,\eps})}
$$
Subsequently, by use of the averaging principle and the Krasnosel'skii type theorem, one can determine 
the fixed point index of the translation operator $\Phi_{T}^{(\eps)}$ for a problem
$$
\dot u(t) = A u (t) + \eps F(t,u(t)), \ t>0.
$$
In particular partial differential equations with a perturbation satisfying some Landesman-Lazer type conditions,
$F$ is the Nemytzkii operator determined by the perturbation, it is possible to employ continuation
along the parameter $\eps\in (0,1]$ and compute the index $\ind (\Phi_T, U)$.
I used such a method for a damped hyperbolic problem in [R5].\\
\indent The already discussed results were applied to study the existence of periodic solutions for partial differential equations. In the paper [R2] we obtained criteria for the existence of periodic solutions for damped hyperbolic equations
without resonance at infinity; the resonant case was examined in [R5]; in the paper [R6] 
I dealt with the beam equation with strong damping; an application to nonnegative periodic solutions of parabolic equations was given in [R4].\\
\indent Periodic problems associated with the equation $(Z)$ or particular partial differential equation are 
subject of intensive studies led by numerous mathematicians and by means of various methods.
Beside the mentioned papers, which are based on translation along trajectories method, it is appropriate to 
refer to such authors as Mawhin \cite{Mawhin},  Fu\v{c}ik and Mawhin \cite{Fucik-Mawhin}, Pr\"{u}ss \cite{Pruss}, Brezis and Nirenberg \cite{Brezis-Nirenberg},
Cesari and Kannan \cite{Cesari-Kannan}, Amman and Zehnder \cite{Amman-Zehnder}, Becker \cite{Becker}, Vrabie \cite{Vrabie}, Hirano \cite{Hirano} and \cite{Hirano1},  Hu and Papageorgiou \cite{Hu-Papa}, Shioji \cite{Shioji1997},
Kamenskii, Obukhovskii and Zecca (\cite{KOZ} and the references therein), Ortega \cite{Ortega}, Mawhin, Ortega and Robles-P\'{e}rez \cite{Mawhin-Ortega-Robles}, Ortega and Robles-P\'{e}rez \cite{Ortega-Robles}.
A short discussion on these other techniques applied there is provided in Section 8 (page \pageref{15062011-1024}).\\
\indent The relation between the fixed point index of the translation operator and 
the topological degree, which was obtained in the dissertation, was also used to prove 
an infinite dimensional version of {\em the Poincar\'{e}-Hopf theorem}.
The classical theorem (see \cite{Milnor} or \cite{Hirsch}) states that if $f$ is a smooth tangent vector field on a compact smooth manifold $M$ and has a finite number of zeros, then the sum of the topological indices of these zeros
is equal to the Euler characteristic $\chi(M)$ of the manifold $M$. 
In other words
$$
\deg (f, M) = \chi(M)
$$
where $\deg(f, M)$ is the topological degree of the vector field 
$f$ with respect to the whole manifold $M$. 
This result has many generalizations where subsets $M$ of a Banach space $\mathbb{E}$,
which are  not differential manifolds and $f:M\to \mathbb{E}$ is a continuous mapping taking values in the tangent cone
(\cite{Clarke}, \cite{Cornet},
\cite{Ben-El-Kryszewski}). Nevertheless, one may ask what if 
$F:\mathbb{E}\to \mathbb{E}$ is not a tangent vector field, that some is trajectories of the equation
$\dot u (t) = f(u(t))$ may leave the set $M$. A general theorem in this direction was achieved by Srzednicki \cite[Th. 5.1]{Srzednicki-formula}. It says that if $f:\R^N\to \R^N$ is a $C^1$ vector field, $B$ is an isolating block
for the flow generated by the equation $\dot u(t)= f(u(t)), \ t\in\R$,  such that 
$B = \overline{ \wew B }$ (i.e. $B$ is the closure of an open set) and both $B$ and the set  of exit points $B^-$ are Euclidean neighborhood retracts, then the Brouwer topological degree $\deg_B ( -f, \wew B)$ 
is equal to the subtraction of the Euler characteristics  $\chi(B)-\chi(B^-).$
At this point a natural question arises. It concerns the validity of this type relation for equations in Banach spaces. An infinite dimensional version of the Poincar\'{e}-Hopf theorem has been recently obtained by Bartsch and Dancer in \cite{Dancer-Bartsch}  for compact vector fields  $I+F$ in a Banach space. In the paper [R3] we provide a first result for general evolution systems in the form (\ref{15042011-1914}). We assume that $A$ generates a compact $C_0$ semigroup and that $F$ is locally Lipschitz and has sublinear growth. There appeared significant difficulties of topological nature that come from the infinite dimensionality of the phase space.

\vspace{15mm}

{\bf 2. Compactness of evolution equations}

We may find some important models for evolution equations in the form $(Z)$ among various types of equations and systems of partial differential equations. As it is well known, topological methods usually need some kind of compactness.
In the case of ordinary differential equations, the compactness of proper mappings is natural due 
to the finite dimension of the phase space. However in the general situation -- for evolution equations in a Banach space -- one imposes proper assumptions on the semigroup generated by the operator $A$ (such as its compactness, equicontinuity or contractivity) and the nonlinear perturbation $F$. Below I discuss some concrete types of equations and their abstract settings. Each of these equations requires a different homotopy invariant and, in consequence, different proof techniques. 

\noindent {\bf Parabolic equations}\\
The first type of compactness can be found in parabolic problems
\be\label{14042011-1144}
\left\{
\begin{array}{ll}
u_t (x,t) = \Delta u(x,t) + f(t,x,u(x,t)), & x\in \Omega, t>0, \\
u(x,t)= 0, & x\in \partial \Omega, t\geq 0,
\end{array} \right.\ee
where $\Omega \subset \R^N$ is a bounded domain with the smooth boundary
$\partial \Omega$, $f:[0,+\infty)\times
\Omega\times \R \to \R$ is continuous and Lipschitz in the third variable
(it assures the existence of solutions for initial value problems).
Let an operator $A:D(A)\to \mathbb{E}$ in $\mathbb{E}:=L^2 (\Omega)$ be given by the formula 
$A\bar u: = \Delta \bar u$, $\bar u\in D(A):=H^2(\Omega)\cap H_{0}^{1}(\Omega)$
(\footnote{$H^k (\Omega):= W^{k,2}(\Omega)$, $k\geq 1$, where $W^{k,2}(\Omega)$ is the Sobolev space , $H_{0}^{1}(\Omega)$ is the closure in the space $H^1 (\Omega)$ of the linear subspace of 
smooth functions with compact support in $\Omega$.}), and  let $F:[0,+\infty)\times \mathbb{E}\to \mathbb{E}$ 
be the Nemytzkii operator determined by $f$, i.e.
$F(t,\bar u)(x):= f(t,x,\bar u(x))$ for a.a. $x\in \Omega$, all
$t\geq 0$ and $\bar u\in \mathbb{E}$.  
With $A$ and $F$ defined in this way, the problem (\ref{14042011-1144}) can be rewritten -- in the abstract form -- 
as the equation $(Z)$.
It follows from the general existence and uniqueness theory for semilinear equations that the operators of translation along trajectories $\Phi_t$, $t\geq 0$, are well defined. 
Observe that, in view of (\ref{15052011-2004}), one has
\be\label{14042011-1125} 
\Phi_t (\bar u) = S_A (t) \bar u + \int_{0}^{t} S_A (t-\tau) F(\tau
,\Phi_\tau (\bar u)) \d \tau \mbox{ for } \bar u\in \mathbb{E}, \, t>0.
\ee
It is known that the semigroup $\{S_A(t):\mathbb{E}\to \mathbb{E}\}_{t\geq 0}$, 
generated by the Laplacian on bounded domains, with the Dirichlet boundary conditions, is {\em compact}
(i.e. $S_A(t)$ is a compact linear operator for all $t>0$),
furthermore, the perturbation $F$ is bounded on bounded subsets.
Consequently, both terms of the right hand side of the formula 
(\ref{14042011-1125}) determine completely continuous mappings. The first one
is compact by definition, whereas the second one can be split into the sum 
\begin{eqnarray*}
\int_{0}^{t} S_A (t-\tau) F(\tau ,\Phi_\tau (\bar u)) \d \tau   & = &
S_{A} (\eps) \left(\int_{0}^{t-\eps} S_A (t-\eps-\tau) F(\tau ,\Phi_\tau (\bar u)) \d \tau\right)\\
& & + \int_{0}^{\eps} S_A (t-\tau) F(\tau ,\Phi_\tau (\bar u)) \d \tau
\end{eqnarray*}
for any $\eps>0$, which implies the compactness of the mapping determined by the second term in (\ref{14042011-1125}). 
Hence, the operators $\Phi_t$, $t>0$, are completely continuous. 
It allows applying the Leray-Schauder fixed point index or, for problems with a constraint set in the phase space $\mathbb{E}$, the fixed point index for compact maps of absolute neighborhood retracts (see \cite{Dugundji-Granas}
and the references therein). The equation $(Z)$ under the above assumptions was investigated in [R3] and [R4].

\noindent {\bf Damped hyperbolic equations}\\
Another type of compactness comes up in second order hyperbolic equations with damping
\be\label{14042011-1134}
\left\{
\begin{array}{ll}
u_{tt} (x,t) + \beta u_t(x,t) = \Delta u(x,t) + f(t,x,u(x,t)),  & x\in \Omega, \ t>0,\\
u(x,t) = 0, & x\in \partial \Omega, t\geq 0,
\end{array} \right.
\ee
where $\Omega$ and $f$ satisfy similar assumptions as before and $\beta>0$ is a damping coefficient.
The problem can be written as an abstract second order equation 
\be\label{15042011-1203}
\ddot u(t) + \beta \dot u(t) + A u(t) + F(t,u(t)) = 0, \ t>0,
\ee
with the operator $A:D(A)\to X$ in the space $X:=L^2(\Omega)$ and the mapping $F:[0,+\infty)\times X\to X$ 
given in a way similar to that before.
The operator $A$ is known to be strictly positive and self-adjoint.
In a space  ${\bold E}:=X^{1/2}\times X^{0}$ we define an operator
${\bold A}:D({\bold A})\to {\bold E}$ by 
\be \label{20042011-1020}
{\bold A} (\bar u, \bar v) = (\bar v, - A \bar u - \beta \bar v), \ \mbox{ for }\ (\bar u, \bar v)\in D({\bold A}):= X^{1} \times X^{1/2},
\ee
where $X^{\theta}$, $\theta\in\R$, stands for the fractional power space determined by $A$.
Further, define a mapping  ${\bold F}:[0,+\infty)\times {\bold E} \to {\bold E}$ by 
\be\label{20042011-1021}
{\bold F} (t, \bar u, \bar v):= (0, F(t,\bar u)),\ \mbox{ for }\ t\geq 0 \mbox{ and } (\bar u, \bar v)\in {\bold E}.
\ee
Then the problem (\ref{15042011-1203}) is transformed into the following equation
\be\label{20042011-1022}
(\dot u(t), \dot v(t)) = {\bold A} (u(t), v(t)) + {\bold F} (t, u(t), v(t)),\ t>0.
\ee
The mapping ${\bold F}$ is completely continuous, due to the compactness of the embedding $X^{1/2}\subset X$,
which in turn is a consequence of the compactness of the semigroup $\{S_A(t):X\to X\}_{t\geq 0}$.
Endowing ${\bold E}$ with a proper scalar product, the operator ${\bold A}$ is strongly $m$-dissipative,
which, according to the Lumer-Philips theorem, stands that ${\bf A}$ generates a $C_0$ semigroup $\{S_{\bold A}(t)\}_{t\geq 0}$ of contractions, i.e. there exists $\omega>0$ such that $\| S_{\bold A} (t) \| \leq e^{-\omega t}$ for all $t\geq 0$.
Therefore the first term in the right hand side of the Duhamel formula analogical to (\ref{14042011-1125}) 
determines a contraction, while the second one defines a mapping that is completely continuous, since 
${\bold F}$ is so. This means that the operators $\Phi_t$, $t>0$, are {\em $k$-set contractions} (or {\em condensing}), i.e.
\be\label{29042011-1246}
\gamma (\Phi_t(V)) \leq  e^{-\omega t} \gamma (V) \mbox{ for any bounded } V\subset {\bold E},
\ee
where $\gamma$ stands for the Hausdorff (or Kuratowski) measure of noncompactness.
For that reason one may apply the Sadovskii version of fixed point index (see \cite{Sadovskii} as well as \cite{Akhmerov-et-al}). This type problems are studied in [R2] and [R5].

\noindent {\bf Beam equation with strong damping}\\
Yet another type of compactness occurs in the beam equation with strong damping
(see (\ref{30062010-1738})), which can be transformed into an abstract form as 
$$
\ddot u (t) + \alpha A\dot u(t) + \beta \dot u(t) + A u(t) +
(a|u(t)|_{1/4}^{2}\! + b + \sigma (A^{1/2} u(t), \dot u(t))_0 ) A^{1/2} u(t)\! =\! f(\omega t), \, t>0,
$$
where $A$ is a strictly positive self-adjoint operator in a Hilbert space $X$  having compact resolvents $(\lma I-A)^{-1}:X\to X$, $\lma>0$, and $f:[0,+\infty]\to X$ is a continuous function. 
Let $(\cdot, \cdot)_{\theta}$ and $|\cdot|_{\theta}$ denote the scalar product and the norm, respectively,
in the fractional power space $X^\theta$ ($\theta\in\R$). 
Define an operator ${\bold A}:D({\bold A})\to {\bold E}$ in a space
${\bold E}:= X^{1/2}\times X^0$ by
$$
\begin{array}{l}
{\bold A} (\bar u, \bar v) := (\bar v, - A (\bar u + \alpha \bar v) - \beta \bar v), \ \mbox{ for } \
(\bar u,\bar v)\in D({\bold A}),\\
D({\bold A}):= \{ (\bar u, \bar v)\in {\bold E} \mid \bar u + \alpha \bar v \in X^1, \, \bar v\in X^{1/2}\},
\end{array}
$$
and a mapping ${\bold F}: [0, +\infty)\times {\bold E} \to {\bold E}$ by the formula
$$
{\bold F}(t,\bar u, \bar v) := f(t)- (a|\bar u|_{1/4}^{2} + b + \sigma (A^{1/2} \bar u, \bar v)_0 ) A^{1/2} \bar u \
\mbox{ for } \ (\bar u,\bar v)\in {\bold E}, \, t\geq 0.
$$
Then the problem can be written in the form (\ref{20042011-1022}) and, similarly as before, 
one can show that, with a properly modified scalar product in ${\bold E}$, ${\bold A}$ generates a $C_0$ semigroup 
$\{S_{\bold A}(t):{\bold E}\to {\bold E} \}_{t\geq 0}$ of contractions.
However in this case the mapping ${\bold F}$ is not completely continuous.
The compactness of the second term in the proper Duhamel formula
analogical to (\ref{14042011-1125}) comes from the fact that 
$S_{\bold A} (t)_{| \{0\}\times X^{0}}$, $t>0$, are compact operators and ${\bold F} ([0,+\infty)\times {\bold E}) \subset \{0\}\times X^{0}$. 
Hence, the operators $\Phi_t:{\bold E}\to {\bold E}$, $t>0$, are $k$-set contractions with respect to the measure of noncompactness, that is the condition (\ref{29042011-1246}) holds. 
Problems with such a structure are subject of [R6].

\vspace{15mm}

{\bf 3. Krasnosel'skii type formula}

We shall consider an autonomous evolution problem
\be\label{25052011-2023}
\left\{
\begin{array}{l}
\dot u (t) = A u (t) +F(u(t)),\, t>0\\
u(t)\in M, t>0
\end{array}\right.
\ee
in a Banach space $\mathbb{E}$ at the presence of the constraint set $M\subset \mathbb{E}$ (in the phase space) 
where  $A:D(A)\to \mathbb{E}$ is the generator of a  $C_0$ semigroup of bounded linear operators on $\mathbb{E}$ and $F:M\to \mathbb{E}$ is a continuous perturbation.
We shall assume that $M$ is a closed convex cone. In applications, the role of the constraint set is usually played by
the cone of nonnegative functions in a proper functional space.
Clearly, if $M=\mathbb{E}$, then the problem reduces to a~problem with no constraints. 
A standard assumption on the operator $A$ is the invariance of $M$ with respect to the resolvents of $A$.
As for $F$ one often assumes about $F$ that $F(M) \subset M$. Nevertheless, this condition is highly restrictive,
since, e.g. it would mean that, in reaction diffusion equations/systems, only chemical reaction producing reagents would be allowed and that in heat equations one could take into account only heat sources.
The right assumption is the so-called tangency condition, which we pose here and which is naturally satisfied 
by nonlinear reaction terms in equations and systems of partial differential equations.\\
\indent Suppose that $A:D(A)\to \mathbb{E}$ is an $m$-dissipative linear operator (\footnote{A linear operator $A:D(A)\to \mathbb{E}$ is called {\em dissipative} provided that
$\|\bar u - \lma A \bar u \| \geq \| u\|$ for all $\lma>0$ and $\bar u\in D(A)$.
A dissipative operator $A$ is said to be {\em $m$-dissipative} if $(I-A)(D(A))=\mathbb{E}$.}) 
generating a~compact $C_0$ semigroup $\{ S_A (t):\mathbb{E}\to \mathbb{E}\}_{t\geq 0}$ on a Banach space $\mathbb{E}$, $M\subset \mathbb{E}$ is a closed convex cone that is {\em invariant } with respect to the resolvents of the operator $A$, i.e.
\be\label{14042011-1034}
(\lma I-A)^{-1} (M)\subset M, \mbox{ for any } \lma>0.
\ee
Additionally, assume that there exists a locally Lipschitz retraction $r: \mathbb{E}\to M$ such that
for some constant $L>0$ 
\be\label{14042011-1035}
\|r(\bar u) - \bar u \| \leq L d_M (\bar u) \mbox{ for all } \bar u \in \mathbb{E}
\ee
where $d_M$ is the distance function (of points) to the set $M$.
Let $F: M\to \mathbb{E}$ be a locally Lipschitz map having sublinear growth  (\footnote{That is there exists a constant $c>0$ such that $\|F(\bar u)\| \leq c(1+\|\bar u\|)$ for all $\bar u\in M$.}) and being {\em tangent } to the set $M$, i.e.
\be\label{25052011-2147}
F(\bar u) \in T_M  (\bar u):= \overline{\bigcup_{\mu>0} \mu (M-\bar u)} \mbox{ for any } \bar u\in M.
\ee
It turns out that, due to the general viability theory of evolution equations, 
the operators of translation along trajectories  $\Phi_t: M\to M$, $t>0$, associated with the equation (\ref{25052011-2023})
% \be\label{25052011-2023}  \dot u(t) = A u(t) + F(u(t)), \ t>0, \ee
are well defined and completely continuous.\\
\indent Under the above assumptions I proved a Krasnosel'skii type theorem.

\noindent {\bf Theorem 3.1} ([R4, Th. 4.1])\\
% \label{14042011-1342}
{\em If an open bounded subset $U\subset M$ is such that
$A \bar u  + F (\bar u) \neq 0$ for $\bar u \in \partial_M U \cap D(A)$, then there exists $t_0>0$ such that
for $t\in (0,t_0]$ the operator $\Phi_t$ has no fixed points in $\partial_M U$ {\em(}the boundary of $U$ relative to $M${\em)} and
\be\label{11082011-1137}
\Ind_M (\Phi_t, U) = \Deg_M ((A,F), U)
\ee
where $\Ind_M$ stands for the fixed point index for compact mappings in $M$ {\em(}see {\em \cite{Dugundji-Granas}}{\em)}
and $\Deg_M$ denotes the topological degree introduced in {\em \cite{JDE2009-Cw-Kr}}.}

\noindent Theorem 3.1 extends my results of \cite{Cwiszewski-JDE2006} where I assumed that the mapping $F$ must take its values in the cone $M$.  I obtained also versions of the above result under different assumptions: in [R1]
(Theorem 5.1, Theorem 5.3) we reject the compactness of the semigroup and prove a formula similar to (\ref{11082011-1137}) (without constraints, i.e. $M=\mathbb{E}$) assuming that $A$ generates a $C_0$ semigroup of contractions in the separable Banach space $\mathbb{E}$ and $F$ is completely continuous (or a $k$-set contraction with respect to the measure of noncompactness); a counterpart of the formula (\ref{11082011-1137}) for the beam equation is provided 
in the paper [R6] (Theorem 3.1).\\
\indent Up to the author's best knowledge, beside the results of \cite{Cwiszewski-JDE2006}, these are the only infinite dimensional results of Krasnosel'skii type.

\noindent {\bf Remark 3.2} (Topological degree for a pair $(A,F)$)\\
(a) If $A$ is a $m$-dissipative operator with compact resolvent, then 
the topological degree can be defined by 
$$
\Deg((A,F), U) := \deg_{LS} (I-(\lma I-A)^{-1} (\lma I + F), U),
$$
for $\lma>0$, where $\deg_{LS}$ denotes the Leray-Schauder topological degree
for compact vector fields. Let us note that if  $\mathbb{E}=\R^n$ and $A=0$, 
then $\Deg ((A,F), U) =\deg (-\lma F, U) = \deg (-F,U)$ where $\deg$ is the Brouwer degree.\\
(b) In the case when we consider the problem with a constraint set $M\subset \mathbb{E}$
and assume merely the tangency of $F$ (as in Theorem 3.1), the previous formula 
does not make sense and breaks down, since the mapping $F$ is defined on the subset $M$
and can take values out of $M$, that is, in general,  the image of $(\lma I+A)^{-1}
(\lma I + F)$ is not contained in $M$. This means that a simple replacement of  the Leray-Schauder index with the  fixed point index for mappings in $M$ is not a right idea. A correct general construction for such pairs $(A,F)$ 
was provided in the paper \cite{JDE2009-Cw-Kr}, where it is assumed  that $M$ is an ${\cal L}$-{\em retract}, i.e. there exists a retraction $r:B(M,\eta)\to M$ of a neighborhood $B(M, \eta):= \{ \bar u\in M\mid d_M (\bar u) <\eta\}$ with some $\eta>0$ satisfying the condition (\ref{14042011-1035}) on $B(M, \eta)$. That degree is given by 
\be\label{15052011-1203} 
\Deg_M ((A,F),U):= \lim_{\lma \to 0^+}
\Ind_M  (\varphi_\lma, U), \ee
where $\varphi_\lma: M\to M$, $\varphi_\lma (\bar u):= (I-\lma A)^{-1} r(\bar u + \lma F(\bar u))$, $\bar u \in M$
(see also [R4, Section 3.2]).\\
(c) In the case where there is a constant $\omega>0$ such that 
$$
\|S_A(t)\| \leq e^{-\omega t } \mbox{ for all } t>0
$$
and there exists $k\in [0, \omega)$ such that 
$$
\gamma (F(V)) \leq k \gamma (V) \mbox{ for any bounded } V\subset \mathbb{E},
$$
in order to get a suitable topological degree we put 
$$
\Deg ((A,F), U) := \deg_{S} (I-(\lma I-A)^{-1} (\lma I+F), U),
$$
where $\lma>0$ is arbitrary and $\deg_S$ stands for the topological degree 
for $k$-set contractions with respect to the measure of noncompactness (see [R1, Section 2]).

\noindent {\bf The proof of Theorem 3.1} is not just a modification of that for the classical finite dimensional 
theorem and had required a new idea. One of the reasons is that in general solutions of the problem
(\ref{25052011-2023}) with an initial value condition $u(0)=\bar u$, $\bar u\in \mathbb{E}$,
are not differentiable at zero  (it is so even if $A$ generates an analytic $C_0$ semigroup and $F\equiv 0$). 
Problems with showing the compactness properties of some homotopies make another obstacle. In the paper [R4],
I propose a proof in which the problem is reduced to equations involving compact vector fields. 
To this end, one considers a homotopy joining the translation operator $\Phi_{t}^{(\lma)}$
for the equation
\be\label{16052011-1536} \dot u (t) = \lma A u(t) + \lma F(u(t)),
\ t>0, \ee 
with the translation operator $\Psi_{t}^{(\lma)}$ for 
\be\label{16052011-1545} \dot u (t) = - u(t) +\varphi_\lma (u(t)), \
t>0, 
\ee 
with a fixed and sufficiently small $\lma>0$ such that 
\be\label{15052011-1245} 
\Deg_M ((A,F), U) = \ind_M (\varphi_\lma,
U) 
\ee 
(see the definition in (\ref{15052011-1203})). 
We indicate the homotopy by considering a family of equations
$$
\dot u (t) = \widetilde A (\mu) u(t) + \widetilde F (u(t),\mu),\  t>0,
$$
where $\widetilde A (\mu): D(\widetilde A (\mu))\to \mathbb{E}$, $\mu \in [0,1]$, are given by
$$
\begin{array}{c}
\widetilde A (\mu) = \mu \lma A + (1-\mu)I, \, \, \mu\in [0,1],\\
D(\widetilde A (\mu)):= D(A), \mbox{ if  } \mu\in (0,1],\  \mbox{ and }\  D(\widetilde A (0))=\mathbb{E},
\end{array}
$$
and  the mapping $\widetilde F : M \times [0,1]\to \mathbb{E}$ is defined by
$$
\widetilde F (\bar u, \mu):=\mu \lma F (\bar u) + (1-\mu)\, \varphi_\lma (\bar u),\ \bar
u\in M, \ \mu\in [0,1].
$$
As It was mentioned, the main difficulty is to get a proper compactness of the already constructed
homotopy. This is caused by a sort of shift of compactness from the linear part to the nonlinearity. 
More precisely, in the equation (\ref{16052011-1536}) the translation operator is completely continuous
due to the compactness of the semigroup generated by $A$, whereas the translation operator $\Psi_{t}^{(\lma)}$ of the equation (\ref{16052011-1545}) is not completely continuous, although the right hand side of that equation is a completely continuous vector field. Nevertheless one may prove that the mentioned homotopy is a $k$-set contraction
with respect to the measure of noncompactness, which allows one to use
the Sadovskii fixed point index -- see  \cite{Akhmerov-et-al}.\\
\indent Hence, for small $\lma>0$ and
$t>0$, the homotopy invariance of fixed point index yields
\be\label{16052011-1523}
\ind_M (\Phi_{t}^{(\lma)}, U) = \ind_M
(\Psi_{t}^{(\lma)}, U).
\ee
Next, in view of the Krasnosel'skii type formula for compact vector fields on convex cones, 
which I proved in \cite{Cwiszewski-JDE2006}, we get
$$
\ind_M  (\Psi_{t}^{(\lma)}, U) = \ind_M (\varphi_\lma, U) \mbox{ for sufficiently small } t>0.
$$
Hence, rescaling the equation (\ref{16052011-1536}) and using both (\ref{16052011-1523}) and (\ref{15052011-1203}) leads, for sufficiently small $t>0$, to the desired formula
$$
\ind_M (\Phi_{\lma t}, U) = \ind_M (\Phi_{t}^{(\lma)}, U) = \ind_M
(\varphi_\lma, U) = \Deg_M  ((A,F), U).
$$

\vspace{15mm}

{\bf 4. Averaging principle for periodic solutions}

Now let us return to the nonautonomous equation $(Z)$, that is
$$
\dot u(t) = A u(t) + F(t, u(t)), \, t>0, \leqno{(Z)}
$$
where $F:[0,+\infty)\times \mathbb{E}\to \mathbb{E}$ is a continuous mapping that is locally Lipschitz in the second variable, has sublinear growth
and is time $T$-periodic, i.e.  $F(t+T, \bar u) = F(t, \bar u)$ for all $t\geq 0$ and $\bar u \in \mathbb{E}$.
We shall search for $T$-periodic solutions by use of averaging principle:
we place the equation $(Z)$ in a family of the following problems $(Z_\lma)$, $\lma \in (0,1]$, 
$$
\dot u(t) = A u(t) + F(t/\lma, u(t)),\, t>0, \leqno{(Z_\lma)}
$$
and next, having standard assumptions on the averaged problem (i.e. the limit problem as $\lma\to 0^+$),
we infer the existence of periodic solutions for small $\lma>0$.
Subsequently, using {\em a priori } type estimates, we shall obtain the existence of periodic solutions for $\lma=1$,
i.e. the original problem $(Z)$.\\
\indent Start with a general averaging principle. Note that, under the above assumptions, for any
$\bar u\in \mathbb{E}$ and $\lma>0$, the equation $(Z_\lma)$ admits the unique solution $u(\cdot; \bar u, \lma):[0,+\infty)\to \mathbb{E}$ satisfying the initial value condition $u(0)=\bar u$. 
The following averaging principle holds.

\noindent {\bf Theorem 4.1} ([R4, Th. 2.4]) \\ %\label{19042011-1148}
{\em Suppose that 
$$
\begin{array}{c}
\mbox{ for any compact $Q_0\subset \mathbb{E}$ and $t>0$ the set }\\
\ \ \ \ \{u(t; \bar u,\lma)\mid \bar u \in  Q_0, \, \lma >0\}
\mbox{ is relatively compact.}
\end{array}\leqno{(C)}
$$
Then
$$
u(t;\bar v, \lma) \to \widehat u (t; \bar u) \mbox{ as } \bar v \to  \bar u, \lma\to 0^+
$$
uniformly with respect to $t$ from bounded intervals,  where $\widehat u (\cdot; \bar u):[0,+\infty)\to \mathbb{E}$ 
is a solution of
$$\left\{\begin{array}{l}
\dot u(t) = A u(t) + \widehat F( u(t)), \ t>0,\\
u (0) = \bar u
\end{array} \right. \leqno{(\widehat Z)}
$$
with the mapping $\widehat F:\mathbb{E}\to \mathbb{E}$ given by 
$\widehat F(\bar u):=\frac{1}{T} \int_{0}^{T} F(t, \bar u) \, dt$, $\bar u\in \mathbb{E}$.}

\noindent The stated theorem is a bit simplified version of the corresponding result of [R4].
Namely, it is assumed there that $A$ and $F$ depend on a parameter and additionally the problem involves a constraint set $M\subset \mathbb{E}$ being invariant with respect to the resolvents of $A$ and the mapping $F$ is tangent to $M$
(see the assumptions (\ref{14042011-1034}) and (\ref{25052011-2147})).
The continuous and compact dependence on parameter plays a key role for using homotopy invariance
(e.g. in the proof of Theorem 4.4).\\
\indent Assumption $(C)$ requires a few words of comment. It appears relatively weak.
One may prove that condition $(C)$ holds in each of the following cases: if $A$ generates a compact 
$C_0$ semigroup, if $F$ is Lipschitz, if $F$ is a $k$-set contraction.
Hence assumption (C) is satisfied for all the mentioned in Section 2 equations.\\
\indent \label{15062011-1001} A first general infinite dimensional averaging principle is due to Henry \cite{Henry} 
for $A$ generating an analytic semigroup. Theorem 4.1 is an extension of Henry's result to the general case  when $A$ is the generator a $C_0$ semigroup, which does not need to be analytic (therefore a different proof is required).
Moreover, Theorem 4.1 is also a generalization of the results obtained by Couchouron and Kamenskii
\cite{Couchouron-Kamenski} (see also \cite{KOZ}) who used different methods and assumed that $F$ satisfies the $k$-set contractivity condition $\gamma (F([0,T]\times V)) \leq k \gamma  (V)$ for any bounded $V\subset \mathbb{E}$ 
and some constant $k\geq 0$. Hence, the Couchouron-Kamenskii version of averaging principle does not apply
in the case where $A$ generates a compact $C_0$ semigroup and the perturbation $F$ does not satisfy additional 
assumption relating compactness. While Theorem 4.1 can be used, since  condition $(C)$ is satisfied. 

\noindent In {\bf the Proof of Theorem 4.1} given in [R4] one considers arbitrary sequences 
$(\bar v_n)$ in $\mathbb{E}$ and $(\lma_n)$ in $(0,+\infty)$ such that
$\bar v_n\to \bar u$ and $\lma_n\to 0^+$ as $n\to +\infty$.
Condition $(C)$ together with proper compactness criteria for solution operators of semilinear equations allows choosing a subsequence of the sequence $(u(\cdot; \bar v_n, \lma_n))$ convergent to some $\widetilde u:[0,+\infty)\to \mathbb{E}$. Then after passage to the limit in the equality
$$
u(t; \bar v_n, \lma_n) = S_A(t) \bar v_n + \int_{0}^{t} S_A(t-s) F(t/\lma_n, u(s;\bar v_n,\lma_n)) \d s,
$$
one gets 
$$
\widetilde u (t) = S_A (t) \bar u + \int_{0}^{t} S_A (t-s) \widehat F(\widetilde u (s)) \, \d s.
$$
It follows from the uniqueness of solutions for $(\widehat Z\,)$ and the fact that the sequences $(\bar v_n)$ and $(\lma_n)$ are arbitrarily chosen that $\widetilde u = \widehat u (\cdot; \bar u)$. Hence the required equality.

Another issue that I dealt with is the existence of branches of periodic solutions 
for problems $(Z_\lma)$. A point  $\bar u_0 \in \mathbb{E}$ is called a {\em branching point} (or {\em cobifurcation point}) provided that  there exists a sequence $(\lma_n)$ in $(0,+\infty)$ and $\lma_n T$-periodic solutions $u_n:[0,+\infty) \to \mathbb{E}$ of $(Z_{\lma_n})$ such that $\lma_n\to 0^+$ and $u_n(0)\to \bar u_0$ as $n\to +\infty$.
The following direct conclusion from the general averaging principle is a necessary condition for branching.

\noindent {\bf Theorem 4.2} ([R4, Th. 5.1])\\ % \label{19042011-1415}
{\em If $\bar u_0\in \mathbb{E}$ is a branching point, then 
$\bar u_0\in D(A)$ and $A\bar u_0+ \widehat F(\bar u_0)=0$.}

\noindent Hence, one should look for branching points among equilibria of the averaged problem.
Earlier less general version of Theorem 4.2 can be found in  \cite[Th. 5.1]{Cwiszewski-JDE2006} and [R1, Th. 6.1].

Naturally, a question about a {\em sufficient} criterion  for the existence of branching points arises.
As one could expect, the answer can be provided by topological methods.

\noindent {\bf Theorem 4.3}\\
{\em Suppose that either  the semigroup $\{S_A(t)\}_{t\geq 0}$ is compact or the space $\mathbb{E}$ is separable,
\be\label{27052011-1125}
\mbox{ there is $\omega>0$ such that $\|S_A(t)\|_{{\cal L}(\mathbb{E},\mathbb{E})} \leq e^{-\omega t}$, $t\geq 0$,}
\ee
and 
\be\label{27052011-1126}
\mbox{ there exists $k\in [0,\omega)$ such that $\gamma ( F([0,T]\times V) )\leq k \gamma (V)$ for bounded $V\subset \mathbb{E}$.}
\ee
Then there exists a branching point in $U$ provided 
$\Deg ((A,\widehat F), U)\neq 0$.}

\noindent  Theorem 4.3 is an immediate consequence of the following important result on averaging of periodic solutions.

\noindent {\bf Theorem 4.4} ([R1, Prop. 4], [R2, Th. 4.4], [R4, Th. 5.5], cf. [R6, Th. 4.2])\\ 
{\em Let the assumptions of Theorem {\em 4.3} hold and  let 
$\Phi_t^{(\lma)}:\mathbb{E}\to \mathbb{E}$, $\lma>0$, be the operator 
of translation along trajectories for $(Z_\lma)$.
Then there exists $\lma_0>0$ such that, for all $\lma\in (0,\lma_0]$,
$\Phi_{\lma T}^{(\lma)} (\bar u) \neq \bar u$ for $\bar u\in \partial U$ and
$$
\Ind (\Phi_{\lma T}^{(\lma)}, U) = \Deg ((A,\widehat F), U).
$$}
The theorem can be viewed as a topological averaging principle for periodic solutions,
since it shows that the existence of periodic solutions of $(Z_\lma)$ (for small $\lma>0$)
depends on the topological index of equilibria of the averaged equation $(\widehat Z)$.\\
\indent Theorem 4.4 can be found in [R1] under the assumption that $\mathbb{E}$ 
is a separable Banach space, the operator $A$ satisfies the condition (\ref{27052011-1125}), while the perturbation $F$ fulfills (\ref{27052011-1126}).
The case where the operator $A$ is replaced by a family of operators $\{ A(t)\}_{t\geq 0}$ 
generating an evolution system is considered in [R2]. Whereas in [R4] 
I proved a version of that theorem on closed convex cone for the operator $A$ generating a compact $C_0$ semigroup.
Moreover, an analogue for an abstract form of the strongly damped beam equation is obtained in [R6].
Theorem 4.4 is an infinite dimensional extension of the results due to Furi and Pera \cite{Furi-Pera} 
treating vector fields on finite dimensional differential manifolds. A different result can be found in \cite{KOZ}
where Kamenskii, Obukhovskii and  Zecca gave a~condition for the existence of branching points, expressed in terms of the topological degree of some averaged mapping acting in a space of periodic functions with values in the space $\mathbb{E}$.

\noindent In {\bf the proof of Theorem 4.4 } one considers a family of homotopies 
$\Theta_\lma: M\times [0,1]\to M$, with a parameter $\lma>0$,  defined by
$$
\Theta_\lma (\bar u,\mu):=u (\lma T),\, \bar u\in \mathbb{E}, \ \mu\in [0,1],
$$
where $u:[0,+\infty)\to \mathbb{E}$ is the solution of 
$$
\dot u (t) = A u(t) + \mu F(t/\lma, u(t)) + (1-\mu) \widehat F(u(t)), \ t>0
$$
satisfying the initial value condition $u(0)=\bar u$. 
Next, using our compactness assumptions and the general averaging principle from Theorem 4.1, one shows by contradiction that there exists $\lma_0>0$ such that, for $\lma\in (0,\lma_0]$ and $\mu\in [0,1]$,  the mappings $\Theta_\lma(\cdot, \mu)$ have no fixed points in the boundary $\part U$. 
Then an application of the homotopy invariance along parameter $\mu$ for a suitable version of fixed point index yields
\be\label{27052011-1200}
\ind (\Phi_{\lma T}^{(\lma)}, U) = \ind (\widehat \Phi_{\lma T}, U), \mbox{ for } \lma\in (0, \lma_0].
\ee
Further, applying the Krasnosel'skii type formula (that is a proper version of Theorem 3.1) 
to the averaged equation $(\widehat Z)$ and eventually decreasing $\lma_0$, we get for small $\lma>0$
$$
\ind (\widehat \Phi_{\lma T}, U) = \Deg ((A, \widehat F), U).
$$
This equality combined with (\ref{27052011-1200}) gives the assertion.

I employed the idea discussed above in the papers [R4] and [R6].
I my opinion, this approach sheds more light on the general structure and the meaning of averaging, compactness properties and the Krasnosel'skii formula, than in the proofs of the analogical theorems that we obtained in [R1] and [R2]. We construct there similar homotopies, but averaging and compactness arguments are "hidden"
in the proofs of the admissibility of these homotopies, that is of the lack of fixed points in the boundary of $U$.

{\em As an application} of Theorem 4.4 w get a sort of {\em continuation principle} 
that appears important from the point of view of applications.

\noindent {\bf Theorem 4.5}  ([R1, Th. 6.2], [R2, Th. 4.10], [R4,  Th. 5.7])\\ % \label{26042011-1431}
{\em Let $A$ and $F$ satisfy the assumptions of Theorem {\em 4.3} and let 
$U\subset E$ be open bounded and  such that\\
{\em (i)} $A\bar u + \widehat F (\bar u)\neq 0$, for $\bar u\in \partial U\cap D(A)$ and $\Deg((A,F), U)\neq 0$;\\
{\em (ii)} for $\lma \in (0,1)$, the problem $(Z_\lma)$ has no
$\lma T$-periodic solutions starting in $\partial U$.\\
Then problem $(Z)$ possesses a $T$-periodic solution starting in $\overline U$.}

\noindent {\bf The proof of Theorem 4.5} is based on an application of the homotopy 
invariance along  parameter $\lma\in (0,1]$ of fixed point index.  
The use of that property is possible due to assumption (ii), which is a sort of the so-called {\em a
priori} estimate and guarantees that the translation operators $\Phi_{\lma T}^{(\lma)}$ have no fixed point in $\part U$.
While in the case $\lma=1$, either there exists a fixed point in the boundary $\part U$ (and  
then one has the desired periodic solution) or the index $\ind(\Phi_{T}^{(1)}, U)$ is well defined and then, by homotopy invariance and Theorem 4.4, we get
$$
\ind(\Phi_{T}^{(1)}, U) = \ind(\Phi_{\lma T}^{(\lma)}, U) = \Deg((A, \widehat F), U), \mbox{ for sufficiently small } \lma>0.
$$
Finally, in view of assumption (i), the fixed point index is nontrivial, which implies the existence of a $T$-periodic solution for $(Z)$.\\
\indent Conditions (i) and (ii) can be verified, for instance, when the mapping $F$ is differentiable
at infinity and the linearization of the right hand side of the equation $(Z)$ has a~trivial kernel -- see [R2, Th. 4.11] and [R4, Th. 5.8]. A part of the results presented in Section 7 were obtained in that way.

\noindent {\bf Remark 4.6 } Instead of the family of problems $(Z_\lma)$ one can consider equations
$$
\dot u (t) = \lma A u(t) + \lma F(t, u(t)), t>0.
$$
It was done so in [R1] and [R2]. Such problems are easily transformed into $(Z_\lma)$ by time rescaling.

\newpage

$\mbox{ }$

\vspace{15mm}

{\bf 5. Averaging principle for hyperbolic equations at resonance }

Consider now  a general problem of the form
\be\label{20042011-1221}
\ddot u (t) + \beta \dot u(t) + A u(t) + F(t,u(t))=0,\, t>0, 
\ee
where $A:D(A)\to X$ is a self-adjoint operator on a Hilbert space $X$, having compact resolvents and such that $\Ker A \neq \{ 0 \}$ and the operator $A+\alpha I$ is strongly positive for some
$\alpha>0$. Here $\beta>0$ and $F:[0,+\infty)\times X \to X$ is a bounded continuous mapping which is locally Lipschitz
with respect to the second variable and $T$-periodic with respect to the first one.
In such a situation we say that the equation is {\em at resonance} at infinity.
By a~solution of (\ref{20042011-1221}) we mean a solution of the problem 
\be\label{20042011-1221a} (\dot
u(t), \dot v(t)) = {\bold A} (u(t), v(t)) + {\bold F} (t, u(t),
v(t)),\, t>0,
\ee
where ${\bold A}$ and ${\bold F}$ are mappings given by the formulae (\ref{20042011-1020}) and (\ref{20042011-1021}),
acting in the space ${\bold E}:=X^{1/2}\times X^0$.\\
\indent Let a mapping $\bar F: N\to N$, with $N:=\Ker A$, be given by
\be\label{30052011-2350}
\bar F(\bar u):= \frac{1}{T}\int_{0}^{T} PF(t,\bar u) \, dt, \, \bar u\in N,
\ee
where $P:X\to N$ is the orthogonal projection onto $N$ (in view of the Riesz-Schauder spectral theory, $\dim N< +\infty$). Let $\Phi_{T}^{(\eps)}:{\bold E}\to {\bold E}$, with a parameter $\eps>0$, be the operator of translation along trajectories for an equation 
$$
(\dot u(t), \dot v(t)) = {\bold A} (u(t), v(t)) + \eps {\bold F} (t, u(t), v(t)),\, t>0,
$$
I proved a theorem that allows determining the fixed point index of the operator $\Phi_{T}^{(\eps)}$ in the resonant case. 

\noindent {\bf Theorem 5.1} ([R5, Th. 3.1]) \\ % \label{29042011-1402}
{\em Let $U\subset N$ be an open bounded set such that 
$\bar F(\bar u)\neq 0$ for all $\bar u\in \part
U$. Then, for any $R, r>0$, there exists $\eps_0>0$ such that,
for all $\eps \in (0, \eps_0]$,
$$
\Ind(\Phi_{T}^{(\eps)}, U\! \oplus\! B_r \times B_R) = (-1)^{k_-}\deg_B (\bar F, U),
$$
with $B_r := \{ \bar u \in \Ker \widetilde P \mid |\bar u|_{1/2}<r\}$,
$\widetilde P:X^{1/2} \to X$ being the orthogonal projection onto $N$,
$B_R:=\{\bar v\in X^{0} \mid |\bar v|_{0} <R  \}$, $\deg_B$ stands for the Brouwer degree and 
$$
k_-:= \sum_{\lma \in \sigma(A)\cap (-\infty,0)} \dim \Ker (A-\lma I).
$$}

\vspace{-5mm}

\noindent The above formula clearly gives a simple criterion for the existence of periodic solutions.

\noindent {\bf Corollary 5.2} ([R5, Cor. 3.6])\\
{\em If $\deg_B (\bar F, U)\neq 0$, then, for sufficiently small $\eps>0$, the problem 
$$
\ddot u(t) + \beta \dot u(t) + A u(t) +\eps F(t,u(t)) = 0, \ t>0,
$$
admits a $T$-periodic solution.}

\noindent The continuation principle below is a corollary of Theorem 5.1, which is
important from the viewpoint of applications.

\noindent {\bf Theorem 5.3} ([R5, Th. 3.7])\\ % \label{23042011-0106}
{\em Suppose that for a number $R_0>0$ one has\\
{\em (i)} $\deg_B (\bar F, B_N (0,R_0))\neq 0$ where $B_N(0,R_0):=\{\bar u\in N\mid |\bar u|_{1/2}<R_0 \}$;\\
{\em (ii)} for $\eps\in (0,1)$, the problem
$$
\ddot u(t) + \beta \dot u(t) +Au(t) +\eps F(t,u(t)) = 0, \ t>0,
$$
\indent  has no $T$-periodic solutions $u$ with $\| (u(0), \dot u(0)) \|_{\bold E} \geq R_0$.\\
Then the equation {\em (\ref{20042011-1221})} admits a 
$T$-periodic solution.}

\indent The assumptions of the above theorem  can be efficiently verified for damped hyperbolic equations in the case when the nonlinearity $F$ satisfies the so-called Landesman-Lazer conditions (see Section 7).

\noindent To {\bf prove Theorem 5.1} one studies a family of problems of the form
\be\label{27052011-1441}
(\dot u (t), \dot v(t)) = \lma {\bold A} u(t) + \lma \eps (0, G(t/\lma, u(t),\mu) ), \ t>0,
\ee
with parameters $\lma\in(0,1]$ and $\mu\in [0,1]$ and the mapping $G:[0,+\infty) \times X^{1/2} \to X^0$
given by
$$
G(t,\bar u, \mu):= - (1-\mu) F(t,(1-\mu)\bar u + \mu \widetilde P\bar u) + \frac{\mu}{T} \int_{0}^{T}
PF(s, (1-\mu)\bar u + \mu \widetilde P\bar u) ds,
$$
for $t\geq 0$, $\bar u\in X^{1/2}$ and $\mu\in [0,1]$. 
Precisely, one has a family of mappings
$\Theta_{T}^{(\lma, \mu)}:{\bold E}\times [0,1] \to {\bold E}$, $\lma>0$ and $\eps \in [0,1]$, defined by
$$
\Theta_{T}^{(\lma, \mu)} ((\bar u, \bar v), \mu):= (u(T), v(T)),
$$
where $(u,v)$ is a solution of the problem (\ref{27052011-1441}) satisfying 
$(u(0), v(0)) =(\bar u, \bar v)$.
It is proved, by  contradiction, that there exists  $\eps_0>0$ such that, for any $\eps\in (0, \eps_0]$
and $\lma\in (0,1]$,
$$
\Theta_{T}^{(\lma,\eps)} ((\bar u, \bar v), \mu) \neq (\bar u, \bar v) \mbox{ for } (\bar u, \bar v) \in \partial {\bold U}
$$
with ${\bold U}:= U\oplus B_r \times B_R$.
This along with the homotopy invariance of fixed point index, first along parameter $\mu$ and subsequently along $\lma$,
yields
$$
\ind (\Phi_{T}^{(\eps)}, {\bold U}) = \ind (\Theta_{T}^{(1,\eps)}(\cdot, 0), {\bold U} )
= \ind (\Theta_{T}^{(1,\eps)}(\cdot, 1), {\bold U}) = \ind (\Theta_{T}^{(\lma,\eps)}(\cdot, 0), {\bold U} ),
$$
for $\lma\in (0,1]$, which means that
$$
\ind (\Phi_{T}^{(\eps)}, {\bold U})  = \ind (\Psi_{T}^{(\lma,\eps)}, {\bold U}),
$$
with the translation operator $\Psi_{T}^{(\lma, \eps)}:{\bold E} \to {\bold E}$ for the equation
$$
(\dot u(t), \dot v(t)) = \lma {\bold A} (u(t), v(t)) + \lma (0, - \eps \bar F (\widetilde P u(t))), \ t>0.
$$
Further, by rescaling, one gets
\be\label{27052011-1327}
\ind (\Phi_{T}^{(\eps)}, {\bold U})  = \ind (\Psi_{T}^{(\lma,\eps)}, {\bold U})
= \ind(\Psi_{\lma T}^{(1,\eps)}, {\bold U}).
\ee
One can show that mild solutions of the equation
$$
(\dot u(t), \dot v(t)) = {\bold A} (u(t), v(t)) + (0, - \eps \bar F (\widetilde P u(t))), \ t>0,
$$
are also mild solutions of an equation 
$$
( \dot u(t), \dot v(t) ) = ({\bold A}+{\bold J}_{\alpha}) (u(t), v(t))
- {\bold J}_{\alpha} (u(t), v(t))+(0, - \eps \bar F (\widetilde P u(t))), \ t>0,
$$
where ${\bold J}_{\alpha}:{\bold E}\to {\bold E}$ is defined by ${\bold J}_\alpha (\bar u, \bar v) = (0,-\alpha \bar u)$,
$(\bar u, \bar v)\in {\bold E}$, and $\alpha>0$ is such that $A+\alpha I$ is strongly positive.
Hence, in view of the appropriate version of the Krasnosel'skii type theorem, one gets
\be\label{27052011-1328}
\ind(\Psi_{\lma T}^{(1,\eps)}, {\bold U}) = \Deg (({\bold A}+{\bold J}_\alpha, -{\bold J}_{\alpha} + (0,-\bar F \widetilde P)),
{\bold U}) \mbox{ for small $\lma>0$ and $\eps>0$}.
\ee
Next, having the product properties of the degree and the formula for the degree of a compact linear vector field, one 
obtains
$$
\begin{array}{ll}
\Deg (({\bold A}+{\bold J}_\alpha, -{\bold J}_{\alpha} + (0,-\bar F \widetilde P)), {\bold U})
& = \deg_{LS} (({\bold A}+{\bold J}_{\alpha})^{-1}(-{\bold J}_{\alpha} + (0,-\bar F \widetilde P)), {\bold U})\\
& = (-1)^{k_-} \deg_B (\bar F, U).
\end{array}
$$
Finally, combining this with (\ref{27052011-1327}) and (\ref{27052011-1328}) the desired formula for the index at the operator $\Phi_{T}^{(\eps)}$ follows.

An analogue for parabolic equations was proved by P. Kokocki in \cite{Kokocki-1} and \cite{Kokocki-2}.

\vspace{15mm}

\noindent {\bf 6. Infinite dimensional Poincar\'{e}-Hopf theorem}

Let $\pi: [0,+\infty) \times \mathbb{E} \to \mathbb{E}$ be a {\em completely continuous semiflow} on a Banach space $\mathbb{E}$, 
%(\footnote{{\em Półpotokiem} $\pi$ na przestrzeni metrycznej $X$ nazywa się ciągłe odwzorowanie
%$\pi:X\times [0,+\infty)\to X$ spełniające warunki $\pi(x,0)=x$ i $\pi (x,t+s) = \pi (\pi(x,t),s)$ dla dowolnych $x\in X$ i $t,s\geq 0$.}),
that is such that $\pi$ sends bounded sets into bounded ones and 
$\pi_t = \pi(t, \cdot)$ is a completely continuous map for each $t>0$. 
If $B$ is  an {\em isolating block} in the sense of Rybakowski (\cite{Rybakowski-TAMS} or \cite{Rybakowski-book}), then
$B^-$ stands for the {\em exit set} (i.e. set of exit points). (\footnote{A function $\sigma:J\to X$, where $J\subset \R$ is an interval, is called a {\em solution for } $\pi$ provided $\pi (s,\sigma (t)) = \sigma (t+s)$ for all $t\in J$ and $s\geq 0$ such that $s+t\in J$.
A point $x\in \part B$ is called a {\em strict egress} ({\em strict ingress}, {\em
bounce-off}, respectively) for $B$ if for any solution
$\sigma:[-\delta_1, \delta_2]\to X$ for $\pi$ satisfying conditions $\sigma(0)=x$,
$\delta_1\geq 0$ and $\delta_2>0$, the following properties hold:\\
\indent (1) there is $\eps_2\in (0,\delta_2]$ such that, for all $t\in (0,\eps_2)$,
$\sigma (t)\not\in B$ ($\sigma (t)\in
\wew B$, $\sigma(t)\not \in B$, resp.);\\
\indent (2) if $\delta_1>0$, then there is $\eps_1\in
(0,\delta_1)$ such that, for all $t\in (-\eps_1,0)$,
$\sigma(t)\in\wew B$ ($\sigma(t)\not\in B$,
$\sigma(t)\not\in B$, resp.).\\
By $B^e$, $B^i$ and $B^b$ one denotes the set of all strict egress, strict ingress and bounce-off points, respectively. 
A closed set $B\subset X$ is said to be an {\em isolating block } for
$\pi$ provided that $\bd B = B^e \cup B^i \cup B^b$ and
the {\em exit set } $B^-:=B^e\cup B^b$ is closed.}) The following result establishes a relation between 
the fixed point index with respect to an isolating block  
and its homotopy type (that is the Conley type homotopy index).

\noindent {\bf Theorem 6.1} ([R3, Th. 1.2 and Th. 4.1])\\
{\em If  $B$ is an isolating block for $\pi$ and $B$ is a neighborhood retract in $\mathbb{E}$, then 
the quotient space $B/B^-$ is of finite (homological) type and there exists $t_0>0$ such that, for $t\in (0,t_0]$, 
$\pi_t:\mathbb{E}\to \mathbb{E}$ has no fixed point in $\part B$ and the following formula holds
$$
\ind (\pi_t, \wew B) = \chi (B/B^-)-1
$$
where $\chi$ denotes the Euler characteristic.}

\noindent Observe that the theorem actually says that the fixed point index of a semiflow
with respect to the interior of an isolating block is equal to the Euler characteristic of the
Conley index in Rybakowski's version. Recall that the Conley index is equal to the homotopy class of the pointed topological space $(B/B^-,[B^-])$ where $B/B^-$ is the quotient topological space and $[B^-]$ is the point relating to the collapsed set $B^-$.

\noindent {\bf The proof of Theorem 6.1} is based on an idea that allows to shift all considerations to the quotient space  $B/B^-$. Namely, we define maps $\pi_t^{B}:B\to B$ by
$$
\pi_{t}^{B}(x):= \pi(\min\{t,s_B(x)\}, x), \ x\in B,
$$
where  $s_B:B\to [0,+\infty]$ is the exit function for the isolating block $B$, i.e.
$s_B(x):=\sup\{t\geq 0 \mid \pi ([0,t]\times \{ 0 \})\subset B \}$, $x\in B$ (see \cite{Rybakowski-TAMS} or \cite{Rybakowski-book}). Then we have $\bar \pi_t: B/B^-\to B/B^-$, $t>0$, given by
$$
\bar \pi_{t} (q(x)) := q(\pi_{t}^{B}(x)), \mbox{ for } x\in B,
$$
where $q:B\to B/B^-$ is the canonical (quotient) projection. This confines our attention to the space $B/B^-$.\\
\indent Another key issue, being a source of essential obstacles, is that our proof requires employing fixed point index theory for maps of $B/B^-$. In the case when the space $\mathbb{E}$ is finite dimensional, $B/B^-$ is a compact absolute neighborhood retract 
(that is it belongs to the class $ANR$) and one can use the Dold index.
However, in this case where $\mathbb{E}$  is of infinite dimension, it appears that the space $B/B^-$ 
in general fails to be metrizable, i.e. it does not belong to the class $ANR$.
This makes it impossible to use the fixed point index in the Granas version for compact maps of spaces from the class $ANR$.
For that reason, we introduced a class of spaces $ANES$ (approximate neighborhood extensor spaces), 
which includes the quotient space $B/B^-$ as well as $ANR$'s. Next we constructed a fixed 
point index for maps of spaces in the class $ANES$.
It turns out that the index possesses many essential properties, that is commutativity, contraction
and normalization, among others. In view of the constructed fixed point index and Lefschetz number, the index $\Ind_{B/B^-} (\bar \pi_t, B/B^-)$ is equal to the Lefschetz number $\Lambda (\bar \pi_t)$ of the map $\pi_t$, which in turn coincides with
the Euler characteristic $\chi (B/B^-)$, since $\bar \pi_t$ is a compact mapping homotopic to the identity.
On the other hand one may show that for sufficiently small $t>0$, $\bar\pi_t$ is constant on some neighborhood of the
point $q(B^-)$ and that the other fixed points in $B$ make a compact subset of $\wew B$. Application of the additivity and commutativity properties leads to the desired formula.

The following generalization of the Poincar\'{e}-Hopf formula
is an immediate and important implication of Theorem 6.1 and the proper version of Theorem 3.1 (see [R3, Th. 5.1]).

\noindent {\bf Theorem 6.2} ([R3, Th. 1.1 and Th. 5.2])\\ % \label{14042011-1357}
{\em 
Let $A$ be the generator of a compact $C_0$ semigroup on $\mathbb{E}$ and let $F:\mathbb{E}\to \mathbb{E}$
be a locally Lipschitz mapping having sublinear growth. If a neighborhood retract
$B\subset \mathbb{E}$ is an isolating block for the semiflow $\Phi:\mathbb{E}\times [0,+\infty)\to \mathbb{E}$ determined by the equation $\dot u(t) =  Au (t) +F(u(t))$, $t>0$, then 
$$
\Deg ((A,F), U)= \chi(B/B^-)-1.
$$}

\noindent Theorem 6.2 is a finite dimensional generalization of the results due to Rybakowski \cite{Rybakowski-formula} 
and Srzednicki \cite{Srzednicki-formula}. Recently there have been obtained a series of theorems in infinite dimension.
Bartsch and Dancer \cite{Dancer-Bartsch} made an extension for vector fields of the form $I+F$ where
$F$ is a compact operator, that is in the case when the compactness occurs in the perturbation.
A result for so-called ${\cal LS}$-vector fields, i.e. fields of the form $L+F$ where $L$ is a 
strongly indefinite bounded operator and $F$ is a compact perturbation, was provided by Styborski
in \cite{Styborski} and \cite{Styborski-phd}. Moreover some results for vector fields $I+F$, where $F$ is a set-valued mapping, which is condensing with respect to the measure of noncompactness, are given by Gudovich, Kamenskii and Quincampoix \cite{GudoKamenQuin2010}. Up to the author's best knowledge, Theorem 6.2 is the only infinite dimensional result with a general $C_0$ semigroup. Moreover, in the mentioned papers \cite{Dancer-Bartsch}, \cite{Styborski} and \cite{GudoKamenQuin2010} it was the perturbation where compactness was assumed, whereas in Theorem 6.2 it is the semigroup generated by $A$ that is compact.\\
\indent After the submission of [R4], prof. R. Srzednicki proposed another proof of Theorem 6.1. In his approach, by use of proper gluing, the space in which mappings are considered remains in the class of $ANR$'s and 
the fixed point index for compact $ANR$'s is sufficient.

\vspace{15mm}

{\bf 7. Averaging and periodic solutions of partial differential equations}

Now we give examples of applications of the abstract results presented in the previous sections to partial differential equations. We shall discuss criteria for periodic solutions for: damped hyperbolic equations, including equations with resonance; beam equations with strong damping  and parabolic equations.
These equations will be formulated as evolution equations of the form $(Z)$ in proper functional spaces, 
and subsequently the abstract theorems from the previous section are to be used.

\noindent {\bf a. Periodic solutions for damped nonlinear hyperbolic equations}

Consider a differential problem
\be\label{21042011-1415}
\left\{
\begin{array}{ll}
u_{tt} (x,t) + \beta(t) u_t (x,t)- \Delta u(x,t) + f(t,u(x,t)) = 0, & x\in \Omega, \ \ t>0,\\
u(x,t) = 0, & x\in \partial \Omega, t>0,\\
u(x,0)=u(x,T), u_t (x,0) = u_t (x, T). & x\in \Omega,
\end{array} \right.
\ee 
where $\Omega\subset \R^N$ is a bounded domain with the smooth boundary 
$\part \Omega$, $\beta:[0,+\infty)\to (0,+\infty)$
is a $T$-periodic function of $C^1$ class and $f:[0,+\infty)\times
\Omega\times \R\to \R$ is a bounded continuous mapping such that 
$$\begin{array}{l}
\bullet \mbox{ for some } T>0, \ \ f(t+T,s)=f(t,s) \mbox{ for all } t>0, \, s\in\R;\\
\mbox{ } \\
 \bullet \mbox{ there is } L>0 \mbox{ such that } |f(t,s_1)-f(t,s_2)| \leq L|s_1-s_2| \mbox{ for  } t\geq 0, s_1, s_2\in\R.
\end{array}
$$
Under these assumptions we proved a criterion for the existence of periodic solutions.

\noindent {\bf Theorem 7.1} (see [R2, Section 5])\\
{\em If the following limit exists
\be\label{17062011-1019}
f_\infty:= \lim_{|s|\to  +\infty} \frac{f(t, s)}{s},
\ee
and the convergence is uniform with respect to $t\geq 0$ and
$f_\infty$ is not an eigenvalue of the Laplacian with the Dirichlet boundary condition, then
the problem {\em
(\ref{21042011-1415})}, admits a solution.}

\noindent {\bf The proof of Theorem 7.1} is based on applying a version of the continuation method from Theorem 4.5
(in this case [R2, Th. 4.10]). According to the discussion in Section 2, the problem (\ref{21042011-1415}) can be rewritten as an abstract second order problem 
$$\left\{
\begin{array}{l}
\ddot u(t) + \beta(t) \dot u(t) + A u(t) + F(t, u(t)) =0, \, t>0,\\
u(0) = u(T), \dot u (0)  = \dot u (T),
\end{array}
\right.
$$
with the operator $A:D(A)\to X$ in $X:=L^2(\Omega)$ given by
$$
A \bar u:= -\Delta \bar u, \ \bar u \in D(A), \ \  D(A) := H^2 (\Omega) \cap H_{0}^{1}(\Omega), \ \ \
$$
and the mapping $F: [0,+\infty) \times X^{1/2} \to X$ defined as 
$$
F(t, \bar u) (x):= f(t, \bar u(x)) \ \ \mbox{ for  a.a. } x\in \Omega \mbox{ and for  all } t\geq 0 \mbox{ and } \bar u \in X^{1/2}.
$$
Observe that the obtained problem can be transformed into a first order periodic problem in the space 
${\bold E}:=X^{1/2}\times X^{0}$
\be\label{26042011-1230}
\left\{\begin{array}{l}
(\dot u (t), \dot v (t)) = {\bold A}(t) (u(t), v(t)) + {\bold F }(t, u(t), v(t)), \ t>0,\\
(u(0), v(0)) = (u(T), v(T)),\end{array}\right. \ee
where the operators
${\bold A}(t)$ in ${\bold E}$ are defined by 
$$
{\bold A}(t)( \bar u, \bar v):=
(\bar v, - A \bar u - \beta(t) \bar v),  \ (\bar u, \bar v)\in D({\bold A}(t)):=X^{1}\times X^{1/2},
$$
and ${\bold F}:[0,+\infty)\times {\bold E}\to {\bold E}$ is given by 
${\bold F} (t, \bar u, \bar v):= (0, F(t,\bar u)), \ t\geq 0, \,
(\bar u, \bar v)\in {\bold E}$.
In view of the assumption (\ref{17062011-1019}),
${\bold F}$ has a derivative at infinity equal to the linear operator ${\bold F}_\infty:{\bold E}\to {\bold E}$,
${\bold F}_\infty (\bar u, \bar v):= (0, -f_\infty \bar u)$, $(\bar u,\bar v)\in {\bold E}$, i.e.
$$
\lim_{\| (\bar u, \bar v)\|_{\bold E} \to + \infty} \frac{{\bold F}(t,\bar u, \bar v) - {\bold F}_\infty (\bar u, \bar v) }{\|(\bar u, \bar v)\|_{\bold E}} =0 \ \ \mbox{ uniformly relative to $t$.}
$$
According to what was said in Section 2, the semigroup generated by ${\bold A}$ (after a proper renorming in ${\bold E}$) has the property $\|S_{\bold A}(t)\| \leq e^{-\omega t}$, for all  $t\geq 0$ and some fixed $\omega>0$, as for the perturbation ${\bold F}$, it is a completely continuous mapping. Therefore the results from Section 4 can be applied.
By the assumption that $f_\infty\not\in \sigma (-A)$, one has 
$\Ker( \widehat {\bold A} + \widehat {\bold F}) =\{ 0\}$ where
$\widehat {\bold A}$ is the average of the family of operators $\{{\bold A} (t)\}_{t\geq 0}$, i.e.
$\widehat {\bold A} (\bar u, \bar v):= (\bar v, -A\bar u - \widehat \beta \bar v)$, $(\bar u, \bar v)\in D(\widehat {\bold A}):= X^{1}\times X^{1/2}$ where $\widehat \beta:= \frac{1}{T}\int_{0}^{T} \beta(t)\d t$.
This implies that, firstly, for sufficiently large $R>0$,
$$
\Deg ((\widehat {\bold A}, \widehat {\bold F} ), B_{\bold E} (0,R)) = \Deg ((\widehat {\bold A}, {\bold F}_\infty), B_{\bold E}(0,R)) = \pm 1,
$$
that is assumption (i) of the continuation principle is satisfied; secondly, the linearized equation 
$$
(\dot u(t), \dot v(t)) = \lma {\bold A}(t)(u(t), v(t)) + \lma {\bold F}_\infty (u(t),v(t)), \ t>0,
$$
does not possess any nontrivial periodic solutions. This property allows proving the existence of $R_0>0$ 
such that,  for any $\lma\in (0,1]$, the problem
$$
\left\{\begin{array}{l}
(\dot u (t), \dot v (t)) = \lma {\bold A}(t) (u(t), v(t)) + \lma {\bold F }(t, u(t), v(t)), \ t>0,\\
(u(0), v(0)) = (u(T), v(T)),\end{array}\right.
$$
does not possess solutions with $\|(u(0), v(0)) \|_{\bold E} \geq R_0$ (see [R2, Proof of Th. 4.11]). 
This means that the {\em a priori } type estimate from assumption (ii) of the continuation principle holds (cf. Th. 4.5 and Remark 4.6). Hence, applying the mentioned continuation principle yields the existence of an (integral) periodic solution for (\ref{26042011-1230}), that is for (\ref{21042011-1415}).

Now turn to the case with resonance and consider a damped hyperbolic problem
\be\label{15042009-1423}
\left\{\begin{array}{ll}
\!\! u_{tt}(x,t) + \beta u_t(x,t) - \Delta u(x,t)  + \lma_k u(x,t) + f(t,x,u(x,t)) \!=\! 0, & x \!\in\! \Omega,         \, t  \! >\!  0, \\
\!\!u(x,t) \!=\! 0,                                               & x\!\in\! \partial \Omega, \, t\! \geq\! 0,\\
\!\! u(x,0)\!=\!u(x,T), u_t (x,0) = u_t (x, T), & x\!\in\! \Omega,
\end{array} \right. \ee
where $\Omega\subset\R^N$ is a bounded domain with the smooth boundary $\partial\Omega$,
$\beta>0$, $\lma_k$ is  the $k$-th eigenvalue of the Laplace operator with the Dirichlet boundary conditions, and $f:\R\times\Omega\times \R\to\R$ is a bounded continuous function satisfying conditions
$$\begin{array}{l}
\bullet\mbox{ there is  $L\!>\!0$ such that
$|f(t,x,s_1)\!-\!f(t,x,s_2)|\!\leq \!L |s_1\!-\!s_2|$},
\mbox{ for {\small $(t,x)\!\in\!\R\!\times \!\Omega$, $s_1,s_2\!\in\! \R$}};\\
\mbox{ }\\
\bullet\mbox{ for some $T>0$, $f(t+T,x,s)=f(t,x,s)$, for all $(t,x,s)\in\R\times \Omega\times \R$.}
\end{array}
$$
For this kind of problems I found a criterion involving {\em Landesman-Lazer type} conditions.

\noindent {\bf Theorem 7.2} ([R5, Th. 4.1])\\ % \label{15042009-1448}
{\em Suppose that one of the following two conditions is satisfied
\be\label{22042009-0201} \int_{0}^{T}\! \left(\,\int\limits_{\{ \phi > 0 \}} \check f_+ (t,x) \phi(x) \d x
+\!\!\int\limits_{\{ \phi<0 \}} \!\hat f_- (t,x) \phi(x) \d x
\right) \d t  > 0 \mbox{ for all } \phi\in N\setminus \{ 0 \},
\ee
where $N$ is the eigenspace of the Laplace operator with the Dirichlet boundary conditions, 
associated with the eigenvalue $\lma_k$, $\check f_+(t,x\! )\!:= \!\liminf\limits_{s\to
+\infty} f(t,x,s)$ and $\hat f_- (t,x\! ) \!:= \!\limsup\limits_{s\to
-\infty} f(t,x,s)$, or 
\be\label{22042009-0202} \int_{0}^{T}\!
\left(\,\int\limits_{\{ \phi > 0 \}} \hat f_+ (t,x) \phi(x)  \d x
+\!\!\int\limits_{\{ \phi<0 \}} \!\check f_- (t,x) \phi(x) \d x
\right) \d t  < 0 \mbox{ for all } \phi\in N\setminus \{0\}, \ee where
$\hat f_+(t,x) \! :=\! \limsup\limits_{s\to +\infty} f(t,x,s)$ and
$\check f_- (t,x)\! :=\! \liminf\limits_{s\to -\infty} f(t,x,s)$.
Then the problem {\em (\ref{15042009-1423})} has a solution.}

\noindent In order to {\bf prove Theorem 7.2}, I applied the continuation principle from Theorem 5.3
and used the Landesman-Lazer type conditions to derive {\em a priori} type estimates.
First, one rewrites the problem  (\ref{15042009-1423}) in the form
$$
\ddot u(t) + \beta \dot u(t) + A u(t) + F(t,u(t)) =0, \ t>0
$$
with the operator $A:D(A)\to X$ in $X:=L^2(\Omega)$ given by
$$
A \bar u:= -\Delta \bar u + \lma_k \bar u, \ \ \bar u\in D(A):= H^2 (\Omega)\times H_{0}^{1} (\Omega)
$$
and $F:[0,+\infty) \times X \to X$ given by $F(t,\bar u)(x):= f(t,x,\bar u(x))$, for a.a. $x\in\Omega$, all 
$t\geq 0$ and $\bar u\in X$. 
By a solution one means a mild solution of the first order problem (\ref{20042011-1221a}). 
The condition (\ref{22042009-0201}) (resp. (\ref{22042009-0202})) enables one to compute the Brouwer topological degree of $\bar F$ given by (\ref{30052011-2350}), that is
one may show that on a sufficiently large ball $B_N (0,R)$, $\bar F$ is homotopic to the identity $I$ (resp. $-I$) 
and then  $\deg_B (\bar F, B_N(0,R))=1$ (resp. $\deg_B (\bar F, B_N(0,R))=(-1)^{\dim\, N}$), 
i.e. condition (i) of the mentioned continuation principle holds.
Furthermore,  by use of the Landesman-Lazer type conditions, one finds  $R_0>0$, such that, for $\eps\in (0,1)$, 
the equation 
$$
\ddot u(t) + \beta \dot u(t) + A u(t) + \eps F(t,u(t)) =0, \ t>0,
$$
does not possess periodic solutions with $\|(u(0), \dot u(0)) \|\geq R_0$, i.e. the condition (ii)
of the continuation principle is satisfied. Accordingly, in view of Theorem 5.3, (\ref{15042009-1423}) admits a~solution.

\indent Theorems 7.1 and  7.2 correspond to the earlier results obtained by different methods
for the telegraph equation (e.g. \cite{Mawhin}, \cite{Fucik-Mawhin}, \cite{Brezis-Nirenberg}).
For a short discussion see Section 8.\\

\noindent {\bf b. Periodic solutions of strongly damped beam equation}

The paper [R6] was motivated by the strongly damped beam equation, that is a problem of the form
\be\label{30062010-1738} \begin{array}{l}
u_{tt}( x,t) \!+\! \alpha u_{txxxx}(x,t) \!+\! \beta u_t(x,t) \!+\! u_{xxxx}( x,t)
 \!- \! \left(a \int_{0}^{l} |u_{\xi}(\xi,t)|^2 d\xi \!+\! b \right) u_{xx}(x,t),  \\
 \ \ \ -  \sigma \left(\int_{0}^{l} u_{\xi}(\xi,t) u_{\xi t}(\xi,t)d\xi \right)  u_{xx}(x,t) = \varphi (x,\omega t),\  \
\!\! \mbox{ for }\! (x,t) \!\in\! (0,l) \!\times\! (0,+\infty),
\end{array} 
\ee 
with the boundary conditions 
\be\label{30062010-1739}
u(0,t)=u(l,t) = 0, u_{xx}(0,t) = u_{xx}(l,t)=0, \mbox{ for } t>0, 
\ee
where $\alpha, \beta, \sigma, a, l >0$, $b\in \R$,
$\varphi:[0,l]\times [0,+\infty)\to\R$ is $T$-periodic with respect to
time variable ($T>0$) and  $\omega>0$ is a parameter. 
The problem (\ref{30062010-1738})--(\ref{30062010-1739}) describes 
-- according to Ball's model (see \cite{Ball}, \cite{Ball1}) -- 
an extensible beam which is subject to an external oscillating force $\varphi$.
It is a special case of the abstract equation
\be\label{07012011-1344} \ddot u + \alpha A  \dot u + \beta \dot u
+ Au + \left(a |u|_{1/4}^{2} + b + \sigma (A^{1/2} u,\dot u)_{0}\right)
A^{1/2}u = f(\omega t )  + \eps f_0 (\omega t), \ee 
where  $A:D(A)\to X$ is a strictly positive self-adjoint operator in a Hilbert space $X$, 
having compact resolvents, $a>0$ and $b\in
\R$ are fixed, $\omega>0$ and $\eps>0$ are parameters,  and $f,
f_0:[0,+\infty)\to X^0$  are  $T$-periodic H\"{o}lder functions.\\
\indent In the next theorem I proved that in the beam equation (\ref{07012011-1344}), under high frequences of the external force, oscillations near equilibrium points appear.

\noindent {\bf Theorem 7.3} ([R2, Th. 5.3, Prop. 5.2])\\ %\label{07012011-1338}
{\em Suppose that
$$
\int_{0}^{T}f(t) \d t = 0  \ \ \ \mbox{ and } \ \ \   |b|\neq \lambda_{j}^{1/2} \mbox{ for all } j\geq 1
$$
and let $k$ be the integer such that $\lambda_{k}^{1/2} < - b < \lambda_{k+1}^{1/2}$ or $k=0$ if $b> - \lambda_{1}^{1/2}$.
Then there exist $\eps_0>0$ and $\omega_0>0$ such that, for all $\eps\in [0, \eps_0]$ and $\omega\geq \omega_0$, 
the equation {\em(\ref{07012011-1344})} has at least $2k+1$ $(T\!/\!\omega)$-periodic solutions.
$u_j^{(\omega, \eps)}$, $j\in \{ i\in \Z \mid |i|\leq k \}$.
Moreover
$$
u_{j}^{(\omega,\eps)} (t) \to \bar u_j  \mbox{ in } X^{1/2} \mbox{ as } \omega\to +\infty, \eps\to 0, \mbox{ uniformly relative to } t\in\R,
$$
and
$$
\dot u_{j}^{(\omega, \eps)}(t) \to 0 \mbox{ in } X^{0} \mbox{ as } \omega\to +\infty, \eps\to 0 \mbox{ uniformly relative to } t\in\R,
$$
for all $j\in \{ i\in \Z \mid |i|\leq k \}$.
Additionally, the topological indices {\em(\footnote{that is the fixed point indices around the starting points of the periodic trajectories for the operator of translation along trajectories by the period of a given solution of the system (\ref{09012011-2248}).})} of the solutions are given by
$$
\Ind ( u_{0}^{(\omega, \eps)} ) = (-1)^{k_0},
$$
where $k_0$ is the integer such that $\lambda_{k_0}^{1/2}<b <\lambda_{k_0+1}^{1/2}$ and $k_0=0$ if $b<\lambda_{1}^{1/2}$,
and
$$
\Ind  ( u_{j}^{(\omega, \eps)} ) = (-1)^{j+1} \ \mbox{ if }\ 0<|j|\leq k.
$$
}

\noindent {\bf To prove Theorem 7.3} one applies a proper version of Theorem 4.4 (see [R6, Th. 4.2]) 
and linearizes the equation at equilibrium points ([R6, Prop. 5.1]).
According to what was said in Section 2, the equation (\ref{07012011-1344}) can be transformed 
into an abstract first order equation in the space ${\bold E}:= X^{1/2}\times X^{0}$
\be\label{09012011-2248} (\dot u(t), \dot v(t)) = {\bold A}
(u(t),v(t)) + {\bold F}(\omega t,u(t), v(t)), \ t>0, 
\ee
where the operator ${\bold A}:D({\bold A})\to {\bold E}$ is defined by 
$$%\be \label{20122010-1641}
\begin{array}{l}
{\bold A} (\bar u, \bar v):= (\bar v, -A (\bar u + \alpha \bar v)
- \beta \bar v), \ \ (\bar u, \bar v)\in D({\bold A}),\\
D({\bold A}):= \{ (\bar u, \bar v)\in {\bold E} \mid \bar u + \alpha \bar v\in X^1, \bar v\in X^{1/2} \},
\end{array}
$$
and the perturbation ${\bold F}:[0,+\infty) \times {\bold E}\times [0,1] \to {\bold E}$ is given by
$$
{\bold F} (t, \bar u, \bar v,\eps)\! :=\! \left(0, \! -\! \left( a|u|_{1/4}^{2} \! +\! b\!+\! \sigma (A^{1/2}\bar u,\bar v)_0\right) A^{1/2}\bar u \!+\!f( t ) \! +\! \eps f_0(t) \right), \ \ t\geq 0, (\bar u,\bar v)\in {\bold E}.
$$
By using the sectoriality of  $-{\bold A}$, one may apply a sort of the so-called one sided estimates 
to prove that (\ref{09012011-2248}) possesses global solutions on $[0,+\infty)$.
Next, employing a proper version of the periodic averaging principle -- [R6, Th. 4.2], with $\lma = 1/\omega$,
along with linearization, yields the assertion.\\

\noindent{\bf c. Nonnegative periodic solutions for parabolic problems }

Consider a parabolic periodic problem
\be\label{28092010-1308}
\left\{ \begin{array}{ll}
u_t (x,t) =  \Delta u(x,t) + f(t,x,u(x,t)), & x\in \Omega, \ t>0, \\
u(x,t)= 0, & x\in \part \Omega,\ t>0,\\
u(x,t)\geq 0, & x\in \Omega,\ t>0,\\
u(x,0)=u(x,T),  & x\in \Omega,
\end{array}
\right.
\ee
where $\Omega \subset \R^N$, $N\geq 1$, is a bounded domain with smooth boundary,  $T>0$ and
$f:[0,+\infty)\times \overline{\Omega} \times [0,+\infty)\to \R$ is a continuous function satisfying the following conditions
\begin{eqnarray}
& \bullet \ & \mbox{there is } L>0 \mbox{ such that } |f(t,x,s_1) - f(t, x, s_2)| \leq L|s_1-s_2| \label{24092010-2240}   \\ \nonumber
& & \mbox{ for } t\geq 0, x\in \overline \Omega \mbox{ and } s_1, s_2\geq 0; \\
& \bullet  &  f(t,x,0) \geq 0 \ \ \   \mbox{ for } t\geq 0, \, x\in \overline \Omega; \label{29112010-2137} \\
& \bullet   & f(t+T, x, s) = f(t,x,s) \ \ \  \mbox{ for } t\geq 0, \, x\in \overline \Omega,\, s\geq 0. \label{24092010-2242}
\end{eqnarray}

Under the above assumptions, I provided a criterion for the existence of nonnegative periodic solutions. 

\noindent {\bf Theorem 7.4}  ([R4, Th. 6.2])\\
{\em Suppose that the conditions  {\em (\ref{24092010-2240})-(\ref{24092010-2242})} hold and that there exists a $T$-periodic continuous  function $f_\infty: [0,+\infty)\to \R$ having the property
\begin{eqnarray}
\label{17102010-2243} \ \ \
& f_\infty(t)=\lim\limits_{s\to +\infty} \displaystyle{ \frac{f(t,x, s)}{s} }\mbox{ uniformly relative to $t\geq 0$, $x\in \Omega$,}\\
& \mbox { } \\
& \displaystyle{\frac{1}{T}} \displaystyle {\int_{0}^{T}} |f_{\infty}(t)| \d t < - \lma_1, \label{27042011-1325}
\end{eqnarray}
where $\lma_1$ is the first eigenvalue of the Laplace operator with the Dirichlet boundary conditions.
Then there exists a nonnegative $T$-periodic solution of the problem {\em (\ref{28092010-1308})}.}

\noindent {\bf The proof of Theorem 7.4} is an application of a proper version of the continuation
principle given in Theorem 4.4 -- [R4, Th. 5.7]. 
First we rewrite the problem  (\ref{28092010-1308}) as an abstract equation
\be\label{28092010-1705}
\left\{\begin{array}{ll}
\dot u(t) = A u(t) + F(t,u(t)), &  t>0,\\
u(t)\in M, & t\geq 0
\end{array}\right.
\ee
in the Hilbert space $\mathbb{E}:=L^2 (\Omega)$, with $M:=\{ \bar u \in \mathbb{E} \mid \bar u \geq 0 \mbox{ a.e. on } \Omega \}$, the operator $A:D(A)\to \mathbb{E}$ given by the formula
$$
A\bar u:= \Delta \bar u, \, \ \ \bar u \in D(A):= H^2 (\Omega)\cap H_{0}^{1}(\Omega),
$$
and the mapping $F:[0,+\infty)\times M\to \mathbb{E}$ defined as
$$
F(t,\bar u)(x):=f(t,x,\bar u(x)) \mbox{ for a.e. } x\in \Omega, \mbox{ for } t\geq 0 \mbox{ and } \bar u\in M.
$$
As it was mentioned in Section 2, $A$ generates a compact $C_0$ semigroup on $\mathbb{E}$.
Furthermore, the cone $M$ is invariant with respect to the resolvents of $A$, i.e.
$(\lma I- A)^{-1}(M)\subset M$ for $\lma>0$ (due to the maximum principle for elliptic operators) -- see [R4, Lem. 6.1 (ii)]. Whereas the assumption (\ref{29112010-2137}) assures that $F$ satisfies the tangency condition,
i.e. $F(t, \bar u)\in T_M (\bar u)$ for all $t\geq 0$ and $\bar u\in M$ ([R4, Lem. 6.1 (iii)]).
$F$ is also locally Lipschitz with respect to the second variable.
Therefore, we can consider the operator of translation along trajectories for (\ref{28092010-1705}).
The assumption (\ref{17102010-2243}) implies the differentiability of $F$ at infinity, that is 
$$
\lim_{\|\bar u\|\to +\infty,\, \bar u \in M} \frac{F(t, \bar u)-F_\infty(t)\bar u}{\|\bar u\|} =0,
$$
with the mapping $F_\infty:[0,+\infty)\to {\cal L} (\mathbb{E}, \mathbb{E})$ defined by 
$$
[F_\infty(t)\bar u](x): = f_\infty (t) \bar u(x), \mbox{ for a.e. } x\in \Omega, \mbox{ for any } t\geq 0, \bar u\in M.
$$
Further, by use of the assumption (\ref{27042011-1325}), one shows that the equation
$$
\dot u (t) = \lma A u(t) + \lma F_\infty (t) u(t), \ t>0,
$$
has no nontrivial $T$-periodic solutions and that 
$\Ker (A+\widehat F_\infty) = \{ 0\}$. 
This enables one to prove that there exists $R_0>0$ such that
$$
\Deg_M ((A, \widehat F),B_M (0,R_0) ) = \Deg_M ((A, \widehat F_\infty), B_M (0,R_0))
$$
(see [R4, Th. 5.8]) and that the problems 
$$
\left\{
\begin{array}{ll}
\dot u(t) = \lma Au(t)+ \lma F(t, u(t)), &  t>0,\\
u(t)\in M & t>0,\\
u(t+T)= u(t) & t>0,
\end{array}
\right.
$$
have no solutions satisfying $\|u(0)\|\geq R_0$.
It follows from the assumption (\ref{27042011-1325}) and the formula for the topological degree with respect to the cone of positive functions (\cite{JDE2009-Cw-Kr}) that
$$
\Deg_M ((A, \widehat F_\infty), B_M (0,R_0))=1.
$$
Hence, employing the proper continuation principle [R4, Th. 5.7], one gets the existence of a $T$-periodic solution
for the problem (\ref{28092010-1705}), that is for (\ref{28092010-1308}).

\noindent {\bf Remark 7.5}\\
The approach presented here, uses translation along trajectories and applies also (without essential modification)
to the case where the Laplace operator is replaced by a general  elliptic operator of second or higher order.
Moreover, the considered boundary conditions can be replaced by others, e.g. Neumann's, Robin's or mixed ones.
It means, in particular, that the general results from previous sections can be applied to the parabolic problem 
\be\label{01062011-2321}
\left\{
\begin{array}{ll}
u_t (x, t) = {\cal A}(x,D) u(x,t) + f(t,x,u(x,t)), &  x\in \Omega,\, t>0,\\
{\cal B}(x,D) u(x,t) = 0, &  x\in \partial \Omega, t>0,\\
u(x,0)=u(x,T), &  x\in\Omega
\end{array}
\right.
\ee
and the hyperbolic one 
\be\label{01062011-2322}
\left\{
\begin{array}{ll}
u_{tt} (x, t) + \beta u_{t} (x,t) = {\cal A}(x,D) u(x,t) + f(t,x,u(x,t)), &  x\in \Omega,\, t>0,\\
{\cal B}(x,D) u(x,t) = 0, &  x\in \partial \Omega, t>0,\\
u(x,0)=u(x,T), &  x\in \Omega,
\end{array}
\right.
\ee
where $f : [0,T] \times \overline \Omega \times \R \to \R$ and  $\beta>0$. 
In the above equations ${\cal A}(x,D)$ stands for a differential operator of the $2m$-th order, $m\geq 1$,
$$
{\cal A}(x,D) = \sum_{|\alpha|\leq 2m} a_{\alpha} (x) D^{\alpha},
$$
where $a_\alpha:\overline \Omega \to \R$, $|\alpha|\leq 2m$, are Lipschitz functions such that the ellipticity condition holds, i.e. there is $\theta>0$ such that 
$$
%(-1)^m 
\sum_{|\alpha|=2m} a_\alpha (x)\xi^{\alpha} \geq \theta |\xi|^{2m} \mbox{ for any } x\in \Omega, \xi\in \R^N.
$$
Whereas ${\cal B}(x,D) = (B_1 (x,D), \ldots, B_m(x,D))$ consists of boundary differential operators 
$B_j$, $j=1, \ldots, m$ in some special form (if $m=1$, then  $B_1$ corresponds to the Neumann, Dirichlet or Robin conditions, if $m>1$ then $B_j = \frac{\partial^{j-1}}{\partial \nu^{j-1}}$, where $\nu$ is the outward normal field on $\partial \Omega$). Then one may define  $ A : D(A) \to \mathbb{E} $ in $ \mathbb{E} := L^2 (\Omega) $ by 
$$
A u:= {\cal A}(x,D)u
$$
with the domain $D(A)$ being the closure in $H^{2m}(\Omega)$ of the set
$$
\{ u\in C^{2m}(\overline\Omega) \mid B_j (x,D) u(x) =0, \, x\in\partial \Omega,\, j=1,\ldots, m \}.
$$
If additionally $A$ is assumed to be self-adjoint, then it is known that $A$ is the generator of a compact (and analytic) $C_0$ semigroup on ${\mathbb{E}} = L^2 (\Omega)$ (see \cite{Tanabe}), which means that 
the discussed translation along trajectories approach can be applied here.

\vspace{15mm}

\noindent {\bf 8. Translation along trajectories method vs. other methods of studying periodic problems}
\label{15062011-1024}

There are a few topological methods that may be used to find periodic solutions of time-dependent partial differential equations. Besides the presented here translation along trajectories approach, we should mention 
methods that use the theory of coincidence degree/index for operators as well as different undetermined coefficients type methods.\\
\indent The problems (\ref{01062011-2321}) and (\ref{01062011-2322}) can be rewritten as a coincidence problem for two operators
\be\label{01062011-2337}
L u = F(u)
\ee
where the operator $L:D(L)\to \mathbb{E}$ is given by
$$
L u = \dot u - A u  \ \ \       \mbox{ or } \ \  Lu = \ddot u + \beta \dot u - A u, \mbox{ respectively},
$$
in a proper space $\mathbb{E}$ of (classes of) functions $u:\Omega\times [0,T]\to \R$,
which satisfy boundary conditions and are $T$-periodic with respect to time, and the nonlinearity 
$F:\mathbb{E}\to \mathbb{E}$
is given by the formula
$$
F(u)(x,t):= f(t,x, u(x,t)), \mbox{ for a.a. } x\in \Omega \mbox{ and } t\in [0,T].
$$
The coincidence problem can be solved by use of a suitable topological degree/coincidence index or variational methods.\\
\indent When the problem (\ref{01062011-2337}) has variational structure, then in order to find solutions one can apply 
critical point theory, e.g. minimizing a suitable functional, mountain pass or linking theorems,
or Morse theory (see e.g. \cite{Amman-Zehnder}). \\
\indent If the problem does not possess variational structure, then one may treat the equation (\ref{01062011-2337}) 
as the fixed point problem $u = L^{-1}F(u)$ provided that the operator $L$ is invertible. If it is not the case, then one adds to both sides of the equation (\ref{01062011-2337}) some perturbation which makes the perturbed operator $L$ invertible (see e.g. \cite{Mawhin}, \cite{Cesari-Kannan}).\\
\indent The second approach, the so-called undetermined coefficients methods, occurs in different variants.
An essential common point of these methods is that the nonlinear part is frozen, that is for any 
$u\in C([0,T],\mathbb{E})$ solutions of the semi-linear problem 
$$
\dot v (t) = A v(t) + f (t), \ t \in [0,T],
$$
with $f:[0,T]\to \mathbb{E}$, $f (t):=F(t, u(t)), t\in [0,T]$ are considered.
In this way, for instance, one may conclude that $T$-periodic solutions of the equation $(Z)$ are fixed points 
of an operator $G:C([0,T], \mathbb{E})\to C([0,T], \mathbb{E})$ given by
$$
G(u)(t):= S_A(t) \left(S_A(T)u(0)+\int_{0}^{T} S_A(T-s) F(s, u(s)) \d s \right) + \int_{0}^{t} S_A(t-s) F(s, u(s))\, \d s
$$
(see e.g. \cite{Vrabie}).
In another variant of this approach, the suitable operator in the space of periodic continuous functions
$C_T([0,T], \mathbb{E})$ is defined as
$$
\widetilde G (u)(t):=S_A(t)(I-S_A(T))^{-1} \int_{0}^{T} S_A(T-s) F(s, u(s))\, \d s + \int_{0}^{t}S_A(t-s)F(s, u(s)) \, \d s.
$$
The latter operator was used in e.g. \cite{KOZ}.

Each of the mentioned methods, i.e. translation along trajectories, coincidence of operators and
undetermined coefficients, has their strong points. The advantage of the two last ones is that they actually do not require uniqueness of solutions for the initial value problem. In the case of the translation operator, the lack of uniqueness implies that its values are no longer points and it is a set-valued mapping. Moreover, the sets being values 
of the translation operator usually are nonconvex. Nevertheless, it is not a~significant obstacle, since 
the fixed point theory of set-valued mappings (\cite{Gorniewicz}, \cite{Kryszewski}) can employed -- see e.g. 
Bothe \cite{Bothe}, Bader and Kryszewski \cite{Bader-Kryszewski} or  \cite{Cwiszewski-JSLN-2006}.
However, the important advantage of translation along trajectories is that it offers a natural view on both described mathematical problems and phenomena that they describe, since considered mappings are in a phase space and 
there are obvious links with dynamical systems. It has quite useful ramifications, e.g. once the index of an operator of translation along trajectories for an isolated periodic solution is determined, we have at least partial information
concerning its stability. It is known that the asymptotic stability of a periodic solution means that the index of the translation operator around the starting point of that solution is equal to one. Hence, either one can exclude the asymptotic stability of that solution or expect the stability and sometimes even deduce it. Results in this direction, i.e. criteria for the asymptotic stability of periodic solutions in terms of the topological index of the translation operator, were shown for the telegraph equation by Ortega \cite{Ortega}.

% Są ciekawe wyniki Ortegi
% \cite{Ortega} dla równania telegrafu, gdzie indeks rozwiązania
% okresowego daje jednoznaczne rozstrzygnięcie stabilności
% rozwiązania (tzn. jeżeli indeks wynosi jeden, to rozwiązanie jes asymptotycznie stabilne).

\newpage

{\bf 9. Discussion of results not included in the dissertation} 
%{Omówienie pozostałych wyników niewchodzących w skład rozprawy}

{\bf Papers published before receiving of the PhD degree:}

\vspace{-4mm}

\noindent [P1] % \bi{16}
\'{C}wiszewski A., Kryszewski W., {\em Equilibria of set-valued maps: a variational approach}, Nonlinear Analysis TMA {\bf 48}, 5 (2002), 707--749.

\noindent [P2] % \bi{15}
\'{C}wiszewski A., Kryszewski W., {\em Approximate Smoothings of Locally Lipschitz  Fun\-ctionals}, Bolletino U.M.I. {\bf 8}, 5-B (2002), 289--320.

\noindent [P3] % \bi{14}
\'{C}wiszewski A., {\em Differential
inclusions with constraints in Banach spaces}, Topological Methods in Nonlinear Analysis, vol. 20 (2002), 119--134.

% \bi{Cwiszewski-doktorat} Ćwiszewski A., {\em Zagadnienia
% różniczkowe z ograniczeniami na stan; stopień topologiczny
% zaburzeń operatorów akretywnych}, rozprawa doktorska, Toruń 2003.

{\bf Papers published after receiving the PhD degree:}

\vspace{-4mm}

\noindent [P4] % \bi{10}
\'{C}wiszewski A., {\em Topological degree methods for perturbations of operators generating compact $C0$ semigroups}, Journal of Differential Equations, Vol. 220 (2) (2006), 434--477.

\noindent [P5] % \bi{11}
\'{C}wiszewski A., Kryszewski W., {\em Homotopy invariants for tangent vector fields on closed sets}, Nonlinear Analysis TMA, vol. 65 (2006), 175--209.

\noindent [P6] % \bi{12}
 \'{C}wiszewski A., Kryszewski W., {\em The constrained degree and fixed-point index theory for set-valued maps}, Nonlinear Analysis TMA, vol. 64 (2006), 2619--2832.

\noindent [P7] % \bi{13}
 \'{C}wiszewski A., {\em Nonlinear Evolution Inclusions with Constraints}, Lecture Notes in Nonlinear Analysis vol. 8 (2006), Juliusz Schauder Center for Nonlinear Studies, 199--238.

\noindent [P8] % \bi{9}
\'{C}wiszewski A., {\em Degree theory for perturbations of m-accretive operators generating compact semigroups with constraints}, Journal of Evolution Equations 7 (2007), 1--33.

\noindent [P9] %\bi{8}
\'{C}wiszewski A., Kryszewski W., {\em Constrained Topological Degree and Positive Solutions of Fully Nonlinear Boundary Value Problems}, Journal of Differential Equations, Vol. 247 (8) (2009), 2235--2269.

\noindent [P10]  % \bi{7}
\'{C}wiszewski A.,  Rybakowski K. P., {\em Singular dynamics of strongly damped beam equation}, Journal of Differential Equations, Vol. 247 (12) (2009),
3202--3233.

{\bf Papers submitted for publication:}

\vspace{-4mm}

\noindent [P11] % \bi{17}
\'{C}wiszewski A., {\em Averaging principle and hyperbolic evolution systems}, 1--20.

{\bf Papers in preparation:}

\vspace{-4mm}

\noindent [P12] % \bi{18}
\'{C}wiszewski A., Kryszewski W., {\em Elliptic partial differential equations with non\-smooth nonlinear parts}.

\noindent [P13] % \bi{19}
\'{C}wiszewski A., Kryszewski W., {\em The fixed point index for c-admissible set-valued maps on
arbitrary absolute neighborhood retracts}.

\noindent [P14] % \bi{20}
\'{C}wiszewski A., Kokocki P., {\em Averaging and
periodic solutions for parabolic equations}.

\noindent [P15] % \bi{21}
\'{C}wiszewski A., Maciejewski M., {\em Positive
solutions for nonlinear boundary value problems with
$p$-Laplacian}.

% \noindent [P16] %\bi{22}
% Ćwiszewski A., Łukasiak R., {\em Averaging for parabolic evolution systems}, w przygotowaniu.

\newpage

My scientific achievements not included in the habilitation concentrate in such areas of nonlinear analysis as:\\[0.4em]
\indent $\bullet\ $ \parbox[t]{120mm}{ evolution equations and inclusions,}\\[0.4em]
\indent $\bullet\ $ \parbox[t]{120mm}{ homotopy invariants for single- and set-valued operators (e.g. 
fixed point index, topological degree, Conley index), including invariants related to problems with state constraints,}\\[0.4em] 
\indent $\bullet\ $ \parbox[t]{120mm}{ nonsmooth analysis and its applications (mainly differential
calculus of locally Lipschitz functions).}

\noindent {\bf Averaging for hyperbolic evolution systems and general averaging principle }\\
In the paper [P11] I gave a general averaging principle for equations of the form
$$
\dot u(t) = A(t/\lma) u(t) + F(t/\lma,u(t)), \ t>0,
$$
where $\{ A(t) \}_{t\geq 0}$ is a family of operators that generates a (general) linear evolution system on a Banach space, and $F$ is a continuous perturbation. It is assumed that the averaging principle in the linear case holds and that
$F$ is almost periodic with respect to the time variable. The obtained theorem is a more general version of Theorem 4.1 of this report and the cited results from \cite{Henry} and \cite{Couchouron-Kamenski}. Moreover, it was proved that the assumption that the linear averaging principle holds is satisfied for the class of hyperbolic evolution systems in the sense of Kato.

\noindent {\bf Singular dynamics of the beam equation}\\
In the paper [P10] we considered a singular passage to the limit in an equation
$$
\eps^2\ddot u(t) +  \eps\delta \dot u(t) +  A (\alpha u(t) + \dot u(t)) + 
\left[ g (|u(t)|_{1/4}^{2}) + \sigma\eps ( u(t), \dot u(t) )_{1/4} \right] A^{1/2} u(t) = 0,
$$
where $A$ is a strictly positive self-adjoint operator with compact resolvents acting in the Hilbert space $X$,
$\alpha>0$, $\sigma, \delta\in\R$ and $g:\R\to \R$. We proved that the dynamics of the above equation is close as $\eps\to 0^+ $ to the dynamics of the equation with $\eps=0$, that is
$$
\dot u(t) = - \alpha A u(t)  - g(|u(t)|_{1/4}^{2})A^{-1/2}u(t),\  t>0.
$$
The singularity entails some difficulties, which come from the fact that when passing to the limit with $\eps\to 0^+$
the type of the equation and even the phase space change essentially. 
We provided some compactness and convergence properties for that singular limit passage.
They together with general results from \cite{Carbinatto-Rybakowski} were used to prove the continuity of the homotopy invariants. Moreover, under an additional dissipativity assumption, we showed the upper semicontinuity
of the global attractors as $\eps\to 0^+$. The results concerning attractors generalize
earlier results of \v{S}ev\v{c}ovi\v{c} \cite{Sevcovic} that where obtained by different methods (based on a specific construction of attractors and strongly using the dissipativity assumption) for the case $\sigma=0$.
We applied different analytic methods that those used in similar papers for damped hyperbolic equations 
\cite{Hale-Raugel} and \cite{Rybakowski-TMNA2003}. Finally, our results give possibility of studying the dynamics of the beam equations without any dissipativity assumption.

\noindent {\bf Equilibria in problems with constraints}\\
The problem of existence of zeros for vector fields on differentiable manifolds can be extended naturally 
to the problem of finding zeros for a continuous $F:M\to \mathbb{E}$ defined on a closed set of constraints
$M\subset \mathbb{E}$ in a Banach space $\mathbb{E}$ under the tangency condition
$$
F(x) \in T_M (x), x\in M,
$$
where $T_M (x)$ is the Bouligand cone of vectors tangent to $M$ at the point $x$. 
In the paper [P1] we studied the case  when $M$ is given in an analytical form, that is 
$M= \{ x\in \mathbb{E} \mid f(x)\leq
0\}$ with a locally Lipschitz functional $f:\mathbb{E}\to \R$.
We introduced classes of {\em regular} and {\em strictly regular} constraint sets and gave criteria sufficient for the existence of zeros for tangent vector fields on $M$ belonging to one of these classes.
Those results along with those of \cite{Ben-El-Kryszewski} significantly generalize the earlier classical ones by
Fan \cite{Fan}, Browder \cite{Browder}, Plaskacz \cite{Plaskacz} as well as those later due to Cornet \cite{Cornet},
Clarke, Ledyaev and Stern \cite{Clarke}, \cite{Clarke-1} and many other authors.
In the paper [P5], we constructed local homotopy invariants detecting zeros of vector fields on constraint sets,
that is we  developed a topological degree theory for vector fields defined on ${\cal L}$-retracts. Recall that the class of ${\cal L}$-retracts, introduced in \cite{Ben-El-Kryszewski}, contains most of sets that can be met in optimization and control theory, that is, among others, closed convex sets, $C^1$ differential manifolds with or without boundary, epi-Lipschitz sets and already mentioned regular sets. The obtained degree is an extension of the degree for vector fields on differentiable manifolds and provides a useful formula that allows determining topological degree of tangent vector fields. In particular, we generalized a classical Poincar\'{e}-Hopf formula stating that if $M$ is a compact
${\cal L}$-retract and a tangent  field $F:M\to \mathbb{E}$ has a finite number of zeros, then
$$
\sum_{x\in F^{-1}(0)} \ind(F,x) = \chi(M),
$$
where $\ind$ stands for the topological index of a zero (in the sense of the already discussed degree), 
and  $\chi(M)$ is the Euler characteristic of the constraint set $M$. 
Furthermore, in the paper [P6] we extended the topological degree onto tangent 
set-valued fields with convex compact values and strongly tangent vector fields with nonconvex values.
As special cases, a series of the known results were derived, in particular, those due to Cornet and Czarnecki.
Moreover, we succeeded to prove that for epi-Lipschitz sets the Euler characteristic can be expressed 
as the topological degree of a set-valued vector field of normal vectors.
In this way we showed that the number introduced by Cornet and used by many other authors is in fact the Euler characteristic (which is crucial for its calculation). We obtained also results in the direction of Deimling's theorem 
\cite{Deimling-inward} (cf. \cite{Deimling}) for vector fields on noncompact convex sets.

\noindent {\bf Non-negative stationary solutions for nonlinear boundary problems and topological degree with constraints}\\
In [P9] we constructed a topological degree detecting solutions of problems having the form 
$$\left\{
\begin{array}{l}
0 \in - A u + F(u)\\
u \in M,
\end{array}\right.
$$
where $M\subset \mathbb{E}$ is a closed set in a Banach space $\mathbb{E}$, $A:D(A)\multimap \mathbb{E}$ is an 
$m$-accretive operator (possibly nonlinear) and $F$ is a tangent continuous perturbation.  
Having that tool,  we find criteria for the existence of solutions for boundary value problems
$$
\left\{
\begin{array}{ll}
-((u'(x))^{p-1})'  = f(x, u(x)),&  x\in (0, L),\\
u(x)\geq 0,  & x\in (0, L), \\
u(0)=(L)=0, &
\end{array}\right.
$$
where $L>0$ and $f:(0,L)\times \R\to\R$ is a continuous mapping 
satisfying some conditions assuring the tangency of the Nemytzkii operator associated with $f$.
These studies are continued in [P15] (also \cite{Maciejewski-doktorat}), where we consider a problem with
the $p$-Laplacian operator
$$
\left\{
\begin{array}{ll}
-div( |\nabla u(x)|^{p-2} \nabla u(x)) = f(x, u(x)), & x\in \Omega\\
u(x)\geq 0, & x\in \Omega,\\
u(x)=0, & x\in \partial \Omega,
\end{array} \right.
$$
where $\Omega\subset \R^N$ ($N\geq 2$) is a bounded domain , $p\geq 2$ and $f$ is a continuous perturbation such that
$f(x, 0) \geq 0$ for a.a. $x\in \Omega$. We find criteria for the existence of nontrivial solutions
expressed in terms of the relation between the principal eigenvalue of the $p$-Laplace operator 
and the asymptotics of $f$ at zero and infinity. We achieve the results, again after having constructed a proper topological degree. The obtained method applies also to systems of such equations without gradient structure (see \cite{Maciejewski-doktorat}).

\noindent {\bf Structure of solution set for initial value problems with constraints\\
 and periodic solutions}\\
In the paper [P8] I provide a topological degree detecting zeros of an operator $-A+F$ in a constraint set
$M$ in a Banach space $\mathbb{E}$, where $A$ is a (nonlinear) $m$-accretive operator generating a compact semigroup of contractions on $\mathbb{E}$ and $F$ is a continuous perturbation. 
The degree is defined as the fixed point index of the operator of translation along trajectories by  sufficiently small time. It was used to study branching of periodic points and continuation principles for constrained problems.
The case without uniqueness property and with the set-valued $F$ is studied in [P7].
The main difficulty in the set-valued (or no uniqueness) setting is that the translation operator generally turns out 
to be set-valued. Hence, in order to apply fixed point index theory one needs a proper regularity of the solution operator. I proved that the solution set of the differential inclusion $\dot u \in - Au + F(u)$ has $R_\delta$ type structure. The results published in [P8] and [P7] were a part of my PhD dissertation, where I also showed similar results with the strongly $m$-accretive $A$ and the perturbation $F$ being a $k$-set contraction with respect to the measure of noncompactness. It was a continuation of the earlier paper [P3], where I studied structure of the set of solutions contained in a constraint set for an inclusion $\dot u\in F(u)$ where $F$ is condensing with respect to the measure of noncompactness.\\
\indent Analogical results to those from [P8] were obtained in [P4] for the operator $A$ being the generator of a compact 
$C_0$ semigroup of bounded linear operators on the Banach space $\mathbb{E}$.
Moreover, I proved there a Krasnosel'skii type theorem. Namely, I showed that, 
if  $M\subset \mathbb{E}$ is a closed convex cone and $F(M) \subset M$, then the 
introduced there topological degree can be expressed as the fixed point index of $(I+\lma A)^{-1}(I+\lma F)$ (for
$\lma>0$). The formula allows effective verification of the assumptions of criteria for the existence of periodic solutions. It was also an inspiration for further research ([R3, [R5]]).
These results generalize Bothe's ones \cite{Bothe}, where only convex constraint sets were allowed, as well as
results due to Bader and Kryszewski \cite{Bader-Kryszewski}, who deal with the linear $A$. In these papers only global homotopy invariants are used.

\noindent {\bf Fixed point index theory for set-valued maps of absolute neighborhood\\ retracts}\\
When studying stationary solutions and periodic problems for nonlinear evolution inclusions in infinite dimensional spaces,
set-valued operators of translation along trajectories may occur.
Namely, there appear mappings $\Phi:M \to 2^M \setminus \{ \emptyset\}$ of
an absolute neighborhood retract $M$ having the form $\Phi = g \circ \Psi$ 
where $g$ is a single-valued continuous mapping
and $\Psi$ is a mapping with values of $R_\delta$ type (i.e. being a {\em cellular set}).
In my PhD dissertation and [P13], the fixed point index for the so-called cellular pairs of maps is constructed.
Any cellular pair is a G\'{o}rniewicz pair and the mentioned translation operator falls into the class of maps represented by cellular pairs. The provided construction fills 
the essential gap in fixed point index theory of set-valued mappings.

\noindent {\bf Approximation of Lipschitz functionals}\\
In many optimization and nonsmooth analysis problems one may encounter
functionals, which are not differentiable but satisfy the local Lipschitz condition.
For such functionals a concept of the so-called Clarke gradient $\partial f$ is considered (it is a notion similar to that of subdifferential for convex functionals). A natural and meaningful from the application viewpoint
is the issue of approximation of locally Lipschitz functionals by means of smooth ones.
In the paper [P2] we proved that a locally Lipschitz functional $f$ can be approximated by $C^1$ functionals
so that the gradient of the approximating functionals is a graph approximation of the Clarke gradient map $\partial f$.
This method was applied in the same paper [P2] to approximate the so-called regular sets of $\R^N$ by means of differentiable manifolds. Whereas in the paper [P12] this method together with variational ones is employed to search for
solutions of elliptic problems with nonsmooth nonlinearities.

\newpage

$\mbox{ }$

\thispagestyle{empty}

\newpage
\noindent {\bf 10. About the author}
\begin{center} 
  Aleksander Kamil \'{C}wiszewski
\end{center}
\noindent {\bf Date \& place of birth:}   May 21, 1975, Szprotawa, Poland.\\[0.5em]
\noindent {\bf Marital status:}                        married, two children.\\[0.5em]
\noindent {\bf Employment}\\[0.5em]
$\bullet$	\parbox[t]{125mm}{since October 2004: adjunct at the Nonlinear Mathematical Aanalysis and 
Topology Chair, Faculty of Mathematics \& Computer Science, UMK in Toru\'{n};}\\[0.5em]
$\bullet$ \parbox[t]{125mm}{October 1999 – September 2004: assistant at the Nonlinear Mathematical Aanalysis and 
Topology Chair, Faculty of Mathematics \& Computer Science, UMK in Toru\'{n};}\\[0.5em]
$\bullet$ \parbox[t]{125mm}{styczeń 1999 – June 1999: assistant in the Mathematical Analysis Department, 
Faculty of Mathematics \& Computer Science, UMK in Toru\'{n}.}\\[0.5em]
\noindent {\bf Education}\\[0.5em]
2003 - \parbox[t]{118mm}{ PhD in Mathematics, UMK in Toru\'{n}, 
            dissertation title: {\em Differential problems with constraints; 
            topological degree for perturbations of $m$-accretive operators}, advisor: prof. dr hab. Wojciech Kryszewski;}\\[0.5em]
1999 – \parbox[t]{118mm}{MSc in Mathematics, UMK in Toru\'{n},
            thesis title: {\em Equilibria of set-valued mappings in constraint sets}, 
            supervisor: prof. dr hab. Wojciech Kryszewski.}

\noindent {\bf Conferences where I had a talk}\\[0.5em]
\noindent $\bullet$  \parbox[t]{125mm}{Symposium on Nonlinear Analysis, Toru\'{n}, Poland, 2011, invited;}\\[0.5em]
$\bullet$ \parbox[t]{125mm}{Summer Conference on General Topology and its Applications, New York, USA, 2011, invited;}\\[0.5em]
$\bullet$ \parbox[t]{125mm}{International Conference on Differential and Difference Equations, 
Ponta Delgada, Portugal, 2011;}\\[0.5em]
$\bullet$ \parbox[t]{125mm}{The 8th American Institute of  Mathematical Sciences Conference on Dynamical
Systems, Differential Equations and Applications, Dresden, Germany, invited;}\\[0.5em]
$\bullet$ \parbox[t]{125mm}{EQUADIFF 12, Brno, The Czech Republic, 2011;}\\[0.5em]
$\bullet$ \parbox[t]{125mm}{World Congress of Nonlinear Analysts, Orlando, USA, 2008, invited;}\\[0.5em]
$\bullet$ \parbox[t]{125mm}{Forum on PDE's, B\c{e}dlewo, Poland, 2008;}\\[0.5em]
$\bullet$ \parbox[t]{125mm}{Symposium on Nonlinear Analysis, Toru\'{n}, Poland, 2007;}\\[0.5em]
$\bullet$ \parbox[t]{125mm}{Workshop on Operator Theory Methods for Differential Equations, Toru\'{n}, Poland, 2006, invited;}\\[0.5em]
$\bullet$ \parbox[t]{125mm}{International Congress of Mathematicians, Madrid, Spain, 2006;}\\[0.5em]
$\bullet$ \parbox[t]{125mm}{International Conference Analysis and Partial Differential Equations, B\c{e}dlewo, Poland, 2006, invited;}\\[0.5em]
$\bullet$ \parbox[t]{125mm}{International Conference on Dynamics, Topology and Computation, B\c{e}dlewo, Poland, 2006;}\\[0.5em]
$\bullet$ \parbox[t]{125mm}{Jim Dugundji Memorial Conference Fixed Point Theory and its Applications, B\c{e}dlewo, Poland, 2005;}\\[0.5em]
$\bullet$ \parbox[t]{125mm}{Seventh International Conference on fixed point theory and applications, Guanajuato, Mexico, 2005;}\\[0.5em]
$\bullet$ \parbox[t]{125mm}{Latest Progress in the Theory of Multifunctions - international conference, B\c{e}dlewo, Poland, 2004;}\\[0.5em]
$\bullet$ \parbox[t]{125mm}{Symposium on Nonlinear Analysis, {\L}\'{o}d\'{z}, Poland, 2004;}\\[0.5em]
$\bullet$ \parbox[t]{125mm}{Topological and Variational Methods in Nonlinear Analysis, B\'{e}dlewo, Poland, 2003;}\\[0.5em]
$\bullet$ \parbox[t]{125mm}{International Conference on Differential and Difference equations, Patras, Greece, 2002;}\\[0.5em]
$\bullet$ \parbox[t]{125mm}{Fixed point theory and applications - International Conference, Haifa, Israel, 2001;}\\[0.5em]
$\bullet$ \parbox[t]{125mm}{Symposium on Nonlinear Analysis, {\L}\'{o}d\'{z}, Poland, 2001;}\\[0.5em]
$\bullet$ \parbox[t]{125mm}{Topological and Variational Methods in Nonlinear Analysis - International
Conference, B\'{e}dlewo, Poland, 2000;}\\[0.5em]
$\bullet$ \parbox[t]{125mm}{Differential and difference equations - International Conf., Voronezh, Russia, 2000;}\\[0.5em]
$\bullet$ \parbox[t]{125mm}{Juliusz Schauder Symposium on Nonlinear Analysis, Toru\'{n}, Poland, 1999.}\\[0.5em] 
\noindent {\bf Other talks}\\[0.3em]
\noindent $\bullet$ \parbox[t]{125mm}{address to The Scientific Committee of J. P. Schauder Center for Nonlinear Studies, Toru\'{n}, Poland, 2008;}\\[0.5em]
\noindent $\bullet$ \parbox[t]{125mm}{address on the occasion of the Scientific Holiday at the Faculty of Mathematics \& Computer Science, UMK in Toru\'{n}, 2007;}\\[0.5em]
\noindent $\bullet$ \parbox[t]{125mm}{seminar talk at The University of Rostock, Germany, 2006;}\\[0.5em]
\noindent $\bullet$ \parbox[t]{125mm}{lecture at the scientific workshops {\em Heat Kernels}, Henrich Fabri Institut, Blaubeuren, Germany, 2006;}\\[0.5em]
\noindent $\bullet$ \parbox[t]{125mm}{seminar talk, Institute of Mathematics, Politechnika Gda\'{n}ska, 2006;}\\[0.5em]
\noindent $\bullet$ \parbox[t]{125mm}{seminar talk, Institute of Mathematics, Uniwersytet Jagiell\'{n}ski, 2003;}\\[0.5em]
\noindent $\bullet$ \parbox[t]{125mm}{Juliusz Schauder lecture series, Uniwersytet Zielonog\'{o}rski, 2002;}\\[0.5em]
\noindent $\bullet$ \parbox[t]{125mm}{seminar talk, Trinity College Dublin, University College Dublin, Ireland, 1999.}\\[0.6em]
\noindent {\bf Longer research stays}\\[0.5mm]
\noindent $\bullet$ \parbox[t]{150mm}{Stay at the University of Rostock (DAAD fellowship), Germany, 2007;}\\[0.5em]
\noindent $\bullet$ \parbox[t]{150mm}{Stay at Trinity College Dublin (TEMPUS fellowship), Ireland, 1999;}\\[0.5em]
\noindent $\bullet$ \parbox[t]{150mm}{Semester of studies at The University of Paderborn, Germany, 1997.}\\[0.5em]
\noindent {\bf Scientific grants \& awards}\\[0.5em]
\noindent $\bullet$ Polish Ministry of Science Research Grant -- investigator, 2010-2012;\\[0.5em]
\noindent $\bullet$ DAAD fellowship grant, 2007;\\[0.5em]
\noindent $\bullet$ Award of the Rector of UMK for scientific achievements, 2007;\\[0.5em]
\noindent $\bullet$ UMK scientific grant - principal investigator, 2007;\\[0.5em]
\noindent $\bullet$  \parbox[t]{125mm}{Polish Science Foundation grant for participation in the International Congress of Mathematicians in Madrid, Spain, 2006;}\\[0.5em]
\noindent $\bullet$ \parbox[t]{125mm}{UMK grant for participation in The School on Nonlinear Differential Equations in Trieste, Italy, 2006;}\\[0.5em]
\noindent $\bullet$  Award of the Rector of UMK for scientific achievements, 2004;\\[0.5em]
\noindent $\bullet$  Polish Ministry of Science Research Grant -- investigator,  2004-2006;\\[0.5em]
\noindent $\bullet$  UMK scientific grant  - principal investigator, 2001;\\[0.5em]
\noindent $\bullet$ \parbox[t]{150mm}{Polish Ministry of Science Grant -- investigator, 2001-2003.}

\end{document}